\documentclass[journal]{IEEEtran}
\IEEEoverridecommandlockouts
\usepackage{graphicx}
\usepackage{multirow}
\usepackage{url}
\usepackage{color}
\usepackage{multirow}
\usepackage{colortbl}
\usepackage[mathscr]{euscript}
\usepackage{amsmath,amssymb,amsbsy}
\usepackage{wrapfig}
\usepackage{fancyhdr}
\usepackage{tikz,pgfplots}

\newcommand\myvspace[1]{}

\newcommand{\CRPS}{{\rm CRPS}}
\newcommand{\cdf}{\mathcal{F}}
\newcommand{\Heavy}{\mathbf{1}}

\ifCLASSINFOpdf

\else

\fi

\newcommand{\dt}{{\Delta t}}
\newcommand{\mat}[1]{\mathbf{{#1}}}
\newcommand{\ipar}{\vec m}

\newcommand{\zz}{\vec{u}}

\newcommand{\mpr} {{\vec{m}}_{\mbox{\tiny pr}} }

\newcommand{\prcov}{\mat{\Gamma}_{\!\mbox{\tiny pr}}}
\newcommand{\ncov}{\mat{\Gamma}_{\!\mbox{\tiny noise}}}
\newcommand{\postcov}{\mat{\Gamma}_{\!\mbox{\tiny post}}}

\newcommand{\like}{\pi_{\!\mbox{\tiny like}}}
\newcommand{\post}{\pi_{\!\mbox{\tiny post}}}
\newcommand{\prior}{\pi_{\!\mbox{\tiny prior}}}

\newcommand{\ipartrue}{\ipar_{\mbox{\tiny true}}}

\newcommand{\obs}{\vec{d}}
\newcommand{\ff}{\vec{f}}                      % parameter-to-obs MAP
\newcommand{\eeta}{\vec{\eta}}
\newcommand{\GM}[2]{\mathcal{N}\!\left( {#1}, {#2}\right)}
\renewcommand{\vec}[1]{{\mathchoice
                     {\mbox{\boldmath$\displaystyle{#1}$}}
                     {\mbox{\boldmath$\textstyle{#1}$}}
                     {\mbox{\boldmath$\scriptstyle{#1}$}}
                     {\mbox{\boldmath$\scriptscriptstyle{#1}$}}}}

\newcommand{\J}{\mathcal{J}}
\newcommand{\R}{\mathbb{R}}

\newcommand{\tidx}{k}

\newcommand{\map} {{m}_{\mbox{\tiny MAP}} }

\newcommand{\cM}[1]{{\color{red} From Mihai: {#1}}}
\newcommand{\cE}[1]{{\color{blue} From Emil: {#1}}}
\newcommand{\cC}[1]{{\color{green} From Cosmin: {#1}}}
\newcommand{\zC}[1]{{\color{brown} From Zheng: {#1}}}
\newcommand{\TTC}[1]{{\color{red} Things to cut: {#1}}}
\definecolor{grassgreen}{RGB}{92,135,39}
\newcommand{\nnote}[1]{\noindent\emph{\textcolor{grassgreen}{N: #1\,}}}

\newcommand{\onetwotable}[2]{#1}

\DeclareMathOperator*{\argmin}{arg\,min}

\graphicspath{{./}}

\begin{document}

\title{A Bayesian Approach for Parameter Estimation with
  Uncertainty for Dynamic Power Systems
%  \\
%  \vspace{-1.4in}
%{\scriptsize \bf Preprint ANL/MCS-P5525-0116}
%\vspace{1.4in}
  %\vspace{-0.4in}
}

\author{No{\'e}mi~Petra,
Cosmin~G.~Petra,
Zheng~Zhang,
Emil~M.~Constantinescu, and
Mihai~Anitescu
\thanks{Cosmin~G.~Petra,
Zheng~Zhang,
Emil~M.~Constantinescu, and
Mihai~Anitescu (e-mail: \{petra,zhengzhang,emconsta,anitescu\}@mcs.anl.gov) are with the Mathematics and
  Computer Science Division, Argonne National Laboratory, Lemont, IL
  60439.}% 
\thanks{Noemi~Petra (e-mail: npetra@ucmerced.edu) is with the
  Department of Applied Mathematics, School of Natural Sciences,
  University of California, Merced, CA 95343.}
\thanks{This material  was based upon work supported in part by the 
Office of Science, U.S. Dept. of Energy, Office of Advanced Scientific
Computing Research, under Contract DE-AC02-06CH11357. N.~P. also acknowledges partial funding
through the U.~S. Dept. of Energy, Office of Workforce Development for
Teachers and Scientists, 2015 Visiting Faculty Program. We thank Hong
Zhang and Shrirang Abhyankar for their help in setting up the IEEE
example and its adjoint implementation.}
\\\myvspace{-0.3in}}

% make the title area
\maketitle

\begin{abstract}
  We address the problem of estimating the uncertainty in the solution
  of power grid inverse problems within the framework of Bayesian
  inference.  We investigate two approaches, an adjoint-based method
  and a stochastic spectral method. These methods are used to estimate
  the maximum a posteriori point of the parameters and their variance,
  which quantifies their uncertainty. Within this framework we
  estimate several parameters of the dynamic power system, such as
  generator inertias, which are not quantifiable in steady-state
  models.  We illustrate the performance of these approaches on a
  9-bus power grid example and analyze the dependence on measurement
  frequency, estimation horizon, perturbation size, and measurement
  noise. We assess the computational efficiency, and discuss the
  expected performance when these methods are applied to large
  systems.
\end{abstract}

\IEEEpeerreviewmaketitle
%
%
%
%

%Inverse problems are the best!

%%%%%%%%%%%%%%%%%
%
\section{Introduction}
%
%%%%%%%%%%%%%%%%%

%\IEEEPARstart{T}he problem of determining the
%\IEEEPARstart{D}etermining
%Determining
Estimating
the
parameters of a
system given noisy measurements is a critical problem in the
operation of energy systems. Decisions about the best and
safe usage of resources depend critically on knowing the
current parameters or states; typically, not all these
quantities are instrumented. Therefore, their values are
obtained indirectly by reconciliation between the mathematical model of
the system and existing measurements by an inverse estimation
procedure, such as state estimation. Before the advent of
phasor measurement units (PMUs)
the phase angle differences in an
electrical network were determined primarily indirectly by 
estimation from SCADA data. While PMU instrumentation can be
rapidly installed on many parts of the power grid,
thus resulting in their phasor angles with respect to a
universal time reference being directly sensed,
the ones without such measurements will still need to be
inferred indirectly from model and measurements.% by using state
%estimation. 

Moreover, the advent of renewable and distributed energy
generation systems creates additional challenges that need
mathematical inversion. The amount,  type, and setting of generation
may not be known a priori by the operator. Therefore, the
parameters of their generator equivalents that need to be
used for balancing the load and assessing the dynamical
stability will need to be determined from measurements. The
dynamical parameters---such as the equivalent inertia of a
windfarm---pose particular challenges because they are not
observable in steady state~\cite{ZhangBankWanEtAl13}. %\cM{Noemi: insert reference to
%  the windfarm estiamtion parameter paper I sent you this
%  summer}
Therefore they are likely to need more frequent data to capture even
fast perturbations from whose transients they can be inverted.  The
rapid deployment of PMU means that such data streams will become
rapidly available, and thus such parameters can be obtained, provided
that dynamic parameter estimation can be carried out.

%In this paper, we investigate inverse problems stemming from
%parameter estimation for energy systems, with a focus on
%dynamical parameter estimation.
In this paper, we focus on dynamical parameter estimation for energy
systems.  Given the increasing dynamic ranges of the energy systems
and the uncertainty due to evolving user behavior and the increased
use of distributed generation, we find it important to provide
uncertainty estimates for these parameters. In this way the operator
can assess the realistic stability range for next-generation energy
systems.

In prior work, parameter estimation in power grid models  typically has been put in the
context of aggregated load models \cite{knyazkin2004parameter}. Most
often the parameters are obtained as a result of least-squares approaches
\cite{choi2006measurement}. Generally, derivative-free methods are
preferred, which typically lead to minimizations based on genetic
algorithms \cite{bai2009novel}; however, derivative-based
least-squares have been introduced by Hiskens et al.~\cite{hiskens1999power,hiskens2001inverse,hiskens2004power}. 

Since, in an operational environment, one needs to provide an answer
in all circumstances, in this work we embrace a Bayesian point of
view. In this case, even with very little information we can produce
an estimate that at least will encapsulate prior information about
the possible ranges of parameters. With more informative data the
estimation will approach the real value of the parameters, without
changing the inference framework.  In this sense the spread of the
posterior probability density function (pdf), namely the solution of
the Bayesian inverse problem, will quantify how much information from
the data can be used for identifying the parameters.  The challenge in
solving this Bayesian inverse problem is in computing statistics of
the pdf, which is a surface in high dimensions. This is extremely
difficult for problems governed by expensive forward models (as is
the power grid model) and high-dimensional parameter spaces (as is the
case for a large-scale power grid). The difficulty stems from the fact
that evaluation of the probability of each point in parameter space
requires solution of the forward problem, and many such evaluations
may be required to adequately sample the posterior density in high
dimensions by conventional Markov-chain Monte Carlo (MCMC)
methods. Hence, quantifying the uncertainties in parameters 
becomes intractable as we increase the grid dimension. Therefore, the
approach we take is based on
%a linearization of the
%parameter-to-observable map, which yields to
a local Gaussian
approximation of the posterior around the maximum a posteriori (MAP)
point. This approximation will be accurate when the
parameter-to-observable map behaves nearly linearly over the support
of the posterior~\cite{Tarantola05}.

%We will
%select as an estimator of the unknown parameters the maximum a
%posteriori (MAP) point. \zC{Need to rewrite the coming sentence to avoid the unfair comparison between MAP and MCMC.} The estimator has the advantage that it can be
%computed using, for example, a gradient-based maximization procedure
%without needing the complexity of Markov-chain Monte Carlo sampling
%methods that is typical of Bayesian exploration of full posteriors.
%%\cM{ to Noemi: I think we should change to MAPE. Then it is a
%%  self-consistent noun. If not OK, please change to MAP point
%%  everywhere}~\nnote{Fixed this.}
%The MAP point can also be used to
%build a Gaussian approximation of the posterior, that will be accurate
%when the parameter-to-observable map behaves nearly linearly over the
%support of the posterior~\cite{Tarantola05}.
%%\cM{to Noemi: discuss
%%  further advantages of MAPE}~\nnote{See if the previous sentence is
% % enough.}
%

We present two  methods for computing MAP and estimating the parametric uncertainty: (1) an adjoint-based
method 
%related to the strategy proposed by Hiskens et
%al. \cite{hiskens2004power} in power grids~\nnote{I don't think we
%  should relate the adjoint-based method to Hisken's paper, from what
%  I see H used symbolic differentiation and sensitivity equations
%  based approach to compute the gradient.} \cE{I agree}
and (2) a surrogate modeling
approach based on polynomial chaos expansions.  These methods solve
the same problem but have different properties and computational
cost. %They are used to compute MAP and to estimate the parametric uncertainty.
We will use these techniques to estimate the inertias
of three generators in an IEEE 9-bus model. The situation
models the circumstance where the actual bulk or distributed inertia is not
known to the grid operator (as would be the case of a
windfarm or other energy resources).

%On
%the other hand, this inertia is necessary in order to
%understand the stability limits of the system following a
%potential contingency, such as a fault-induced transmission
%line relay trip, which is a required reliability analysis
%for all system operators.

We carry out extensive validation experiments to demonstrate the
consistency and accuracy of the methods. Moreover, we use our approach
to investigate the effect of important data features on the precision
of MAP. These features include the frequency of the measurements and
the size of the perturbation.  In addition, we compare the behavior of
the two methods and discuss their computational efficiency and the
expected complexities when applied to larger systems.
%what we can expect for their
%complexity when these methods are
%applied to larger systems. 

%The model used in this work is implemented in PETSc along with
%Jacobian information \cite{petsc-web-page,petsc-user-ref} and
%available as an open-source example.~\nnote{This previous sentence is
%  hanging a bit, perhaps we can move it or remove it.}

%%%%%%%%%%%%%%%%%
%
\section{Problem Formulation}\label{s:probFor}
\begin{figure}
  \centering\includegraphics[width=0.62\columnwidth]{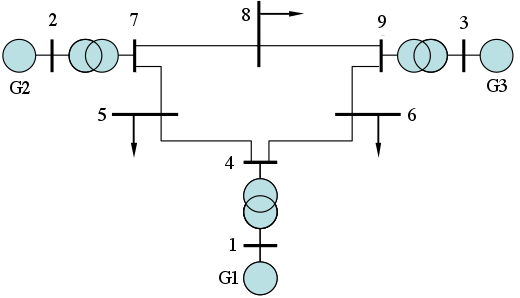}
  \myvspace{-0.1in}
\caption{IEEE 9-bus test case system. Here the buses 1, 2 and 3
  are generator buses and 5, 6, and 8 load buses.}
\myvspace{-0.3in}
\label{fig:3gen9bus}
\end{figure}

%We use the $9$-bus power system as a demonstration
Assume that we have measurements of a dynamical system that can
be modeled by an additive Gaussian noise model
\begin{equation}\label{equ:noise-model}
    \obs = \ff(\ipar) + \vec{\eta}, \quad \eeta \sim
    \GM{\vec{0}}{\ncov},
\end{equation}
where $\ncov\in \R^{q\times q}$ is the measurement noise covariance
matrix and $\ff(\cdot)$ a (generally nonlinear) operator mapping model
parameters $\ipar$ to observations $\obs$.  Here, evaluation of this
{\it parameter-to-observable} map $\ff(\ipar)$ requires solution of a
differential-algebraic system (DAE) that models the dynamics of a
power grid (followed by extraction of the DAE solution at observation
points):
%$\ff$ denotes the observable quantities $\obs$ that depend on parameters $\ipar$.% to measurements $\obs$.
%%\zC{I would suggest saying
%%  ``f denotes some observable quantities dependent on parameters
%%  $\ipar$".
%%It's easier to understand for engineers.}.
%The function $\ff$ is
%computed  by solving the following
\begin{subequations}\label{eq:powersystem}
\begin{alignat}{2}
  \dot{\vec x} &= h(t,\vec x,\vec y,\ipar)\,,~~&&\vec x(0) = \vec x_0\,,~\label{eq:ode}\\
  \vec 0 &= g(t,\vec x,\vec y)\,, ~~  \quad  &&\vec y(0) = \vec y_0\,. \label{eq:algeqs} 
\end{alignat}
\end{subequations}
Here $\vec x$ represents the dynamic state variables (e.g., rotor
angle, generator speed), $\vec y$ the static algebraic variables
(e.g., bus voltages and line currents), $\vec x_0$ and $\vec y_0$ are
the initial states, $t$ represents time, and $\ipar$ the model
parameters.  The right-hand side $h$ in~\eqref{eq:ode} is in general a
nonlinear function that models the dynamics of the system, and $g$ in
\eqref{eq:algeqs} is a set of algebraic equations modeling the passive
network of the power system. For the IEEE 9-bus power grid model
problem, as illustrated in Figure~\ref{fig:3gen9bus}, for each
generator we have seven differential (i.e.,~$\vec x \in \R^{21}$) and
two algebraic equations, and for each network node two additional
algebraic equations (i.e.,~$\vec y \in
\R^{24}$)~\cite{SauerPai98}. The inference parameter $\ipar$ we
consider in this paper is the inertia of each generator, and thus
$\ipar\in\R^{3}$.  In realistic applications, the initial state
$\vec x_0$ may not be known either, and it would have to be inferred
from data.  However, since our focus is on understanding the
reconstructability of parameters that cannot be determined from
steady-state measurements, such as inertias, we assume that the initial
conditions $\vec x_0$ are known. %\cite{hiskens2004power}.
Nevertheless, initial conditions can also be considered uncertain; and
the framework introduced herein naturally extends to such cases.%, by
%redefining the mapping $\ff$.
%\cM{Aye aye. The initial state is not known any more than
%  the parameters, and actually, less. We should have tried
%  to invert it as well from mixed PMU scada data or at least
%  to look at incompletely PMU instrumented setups. We may get grief for this. No need
%  to answer this, but a heads up} \cE{Hiskens does that but in a
%  slightly different context (I left it commented out). We can say for
%sure that inverting for the initial conditions is also fine for the
%tools introduced here -- so I
% added a sentence to that effect.}
%We use $\vec u := (\vec x,\vec y)$ to denote the state variables.

The Bayesian formulation poses the parameter estimation problem as a
problem of statistical inference over parameter space. The solution of
the resulting Bayesian inverse problem is a {\it posterior probability
  density function (pdf)}.
%that represents our belief about the correct value of
%the parameter.
Bayes' Theorem states the posterior pdf explicitly as
%To compute this distribution, we use Bayes' theorem,
%namely
%\vspace{-0.05in}
\begin{align}
  \nonumber %\label{eq:posterior}
\post(\ipar)
\! \propto \like(\obs | \ipar ) \; \prior(\ipar).
\end{align}
%Solving the parameter estimation problem
%using Bayes' approach
%requires, therefore, specification of the {\it likelihood} 
%, which
%quantifies the relative probability that a candidate parameter $\ipar$
%could have produced the observed data $\obs$,
%and {\it prior} models.
%, which
%describes our beliefs about the parameter before any data are
%considered.

For the likelihood model, $\like(\obs | \ipar )$, we assume that the measured quantities are
the bus voltages from a disturbance. We note that here we measure the
voltage at all buses in the IEEE 9-bus power grid; however, our
framework can be used to experiment with various measurement scenarios
(e.g., measurements at a subset of buses) at various time intervals
and measurements of different quantities. %\TTC{I think starting from
%``More concretely", all
 % sentences in this paragraph can be removed. IEEE engineers will not
 % care about the abstract mathematical operators, and they only care
 % about the physical meanings of $\ff$ and $\obs$ which have already
 % been explained in the beginning of this section.}
%~\nnote{Z's suggestion might work.}
%More concretely, we define a {\em network-time observation
%  operator} $\B : \R^{s} \to \R^{q}$ that projects the DAE state
%solution vector onto the observable vector.
%\cM{To Noemi and Cosmin:
%  We need to describe what is measured. Voltages at all buses?  That
%  is what is the operator ${\cal B}$?.}~\nnote{Fixed this. Also moved
%  this paragraph after the discussion on uncertain params to avoid
%  breaking the flow.}
%Therefore, the parameter-to-observable map
%$\ff$ is a nonlinear operator that maps a parameter vector $\ipar \in
%\R^n$ to the network-time observation vector $\obs \in \R^{q}$, namely,
%%\begin{equation}\label{equ:param-to-obs}
%$\vec f:\vec{m}
%   \,
%   \stackrel{\sol}{\longmapsto}
%   \,
%   \uu
%   \,
%   \stackrel{\B}{\longmapsto}
%   \,
%   \obs$,
%%\end{equation}
% where $\sol$ is the DAE discretization operator  and $\uu \in
%\R^{s}$ is the discrete DAE solution vector.
Furthermore, since the noise $\vec{\eta}$ is independent of $\ipar$, thus
$\obs | \ipar \sim \GM{\ff(\ipar)}{\ncov}$, the likelihood is given by
%\begin{equation} \label{eq:likelihood}
%\vspace{-0.05in}
\begin{align} \label{eq:likelihood}
\like(\obs | \ipar) \propto \exp \Big( -\|\vec f(\ipar) -
  \obs) \|^2_{\ncov^{-1}} \Big).
\end{align}
%\end{equation}
%According to Bayes' theorem
%with Gaussian noise and prior,
%the
% posteriori density function of $\ipar$ is described
% as
 %~\cite{Tarantola05,KaipioSomersalo05}
%\begin{align}
%  \label{eq:posterior}
%\post(\ipar)
%\! \propto \! \exp \Big (\!\!-\!\frac{1}{2}\|\vec f(\ipar)\!-\!\obs)
%\|^2_{\ncov^{-1}} \!\!-\!\frac{1}{2}\|\ipar\!-\!\mpr\|^2_{\prcov^{-1}}\!\Big).
 %\end{align}
%Here $\mpr$ and $\prcov\in \mathbb{R}^{n\times n}$ are the mean and
%covariance matrix of the prior distribution $\rho_{\rm prior}(\ipar)$, respectively.
%
 %, and
%$\ncov\in \mathbb{R}^{q\times q}$ is the noise covariance matrix.
%here taken as $\diag(\sigma_1^2,
%\sigma_2^2, \ldots, \sigma_q^2)$.
%
%\cM{Noemi, you will have to check me on the following, and explain
%  how you chose the prior and the noise}
%
%\TTC{I think we need to shorten this paragraph by making the descriptions more brief.} 
%To set up the Bayesian framework, we need to specify the noise covariance
%$\ncov$,  prior covariance $\prcov$, and prior mean $\mpr$.
%The first can be obtained
The noise covariance, $\ncov$, can be obtained by offline studies of the
measurement setup. If the measurements are from PMUs, one can
reasonably assume that the measurement noise is independent between
sensors and white noise in time for one of them (on the time scale of
interest, which is between 0.03 and 30 s). The variance then can be
computed from the precision rating of the instrument.
%This is consistent with our choice of noise covariance
%above.

The Bayesian prior, $\prior(\ipar)$, on the other hand, requires quantification
of the existing information about the parameters. Considerable
literature exists in the area of eliciting priors, but it certainly
requires an intimate analysis of the system at
hand~\cite{OakleyO'hagan07}. %,Bui-Thanh12}.
%\cM{Noemi, please insert
%  bibliography}~\nnote{Added two refs. The second one is only a report
%  but it has a nice section on choosing the prior.}
For example, for a windfarm, one can use historical logs or
a simulation-based model to create a statistical model of the active
inertia at a given time of the year, conditional on ambient
conditions, or use information from similar windfarms.
%The resulting
%distribution can become the prior.
Here we use a Gaussian prior,
a common choice for Bayesian inverse
problems~\cite{KaipioSomersalo05}.
%The prior mean describes our best
%guess about the uncertain parameter, which could be obtained from
%existing measurements or from other available information.
We use a
prior with large variance and diagonal covariance because of the lack of a priori information
about the parameters.
%; and we assume that the parameters are
%uncorrelated, which essentially leads to a diagonal prior covariance
%matrix.
%\cM{Noemi, explain how you came up with the prior an
%  why}~\nnote{Fixed this. Emil \& Cosmin, please see if you would like
%  to add change anything here. In particular if you have a good
%  explenation on why we can assume that the inertia is uncorelated,
%  that would be great. Also, it would be good to make this shorter.}

Restating Bayes' theorem with Gaussian noise and prior, the
posteriori density function of $\ipar$ is described as
~\cite{Tarantola05,KaipioSomersalo05}
\begin{align}
  \label{eq:posterior}
\post(\ipar)
\! \propto \! \exp \Big (\!\!-\!\frac{1}{2}\|\vec f(\ipar)\!-\!\obs)
\|^2_{\ncov^{-1}} \!\!-\!\frac{1}{2}\|\ipar\!-\!\mpr\|^2_{\prcov^{-1}}\!\Big),
\end{align}
where $\mpr$ and $\prcov\in \mathbb{R}^{n\times n}$ are the mean and
covariance matrix of the prior distribution $\prior(\ipar)$,
respectively.

Despite the choice of Gaussian prior and noise
probability distributions, the posterior probability distribution need
not be Gaussian, because of the nonlinearity of $\vec
f(\ipar)$~\cite{Tarantola05,KaipioSomersalo05}.
%%We are thus led to make a quadratic approximation of the negative log
%%of the posterior \eqref{eq:posterior}, which results in a Gaussian
%%approximation of the posterior $\mathcal{N}(\vec{\map},\postcov)$.
%The maximum a-posterior parameter estimate, $\vec{\map}$, is the
%parameter vector maximizing the posterior \eqref{eq:posterior}.  It can
%be found by minimizing the negative log posterior, which amounts to
%solving the following optimization problem
%\begin{equation}\label{eq:objfunction-bayesian}
%\vec{\map} = \underset{\ipar}{\argmin}
%\; \J(\ipar) := - \log \pi_{\text{post}}(\ipar).
%\end{equation}
Here we make a quadratic approximation of the negative log of the
posterior \eqref{eq:posterior}, resulting in a Gaussian
approximation $\mathcal{N}(\vec{\map},\postcov)$ of the posterior.
%In
%the vast majority of cases stemming from physical systems data, this
%is a very reasonable assumption.
%\cM{Noemi, references?}~\nnote{Hm,
%  not sure if we can add a ref for this statement. But perhaps is better to say:
%  ``A Gaussian is often a good approximation to the posterior when a
%  nonlinear parameter-to-observable map is well approximated by a
%  linearization over the set of parameters with significant posterior
%  probability~\cite{PetraMartinStadlerEtAl14}.''}
The mean of this
posterior approximation, $\vec{\map}$, is the MAP point obtained by minimizing
the negative log posterior:
\begin{align}\label{eq:objfunction-bayesian}
\vec{\map} = \underset{\ipar}{\arg \min}
\; \J(\ipar) := - \log \post(\ipar).
\end{align}
The posterior covariance matrix $\postcov$ is then obtained by computing the
inverse of the Hessian of $\J$ at $\ipar$~\cite{Tarantola05,KaipioSomersalo05}.

%\TTC{This paragraph can be removed.}%~\nnote{May be possible.}
%Solving \eqref{eq:objfunction-bayesian} and finding these
%quantities complete the estimation and provide uncertainty
%bars around the parameters of interest, in our
%example, generator inertias. %Solution methods for this problem are presented in the next section. 

%%%%%%%%%%%%%%%%%
%
\section{Solution Methods}\label{s:solmet}

We present two methods for solving the inverse problem: the
adjoint-based method and the stochastic 
spectral method. In Section~\ref{s:numerical} we will illustrate the circumstances in
which one approach will be favored over the other. 

\myvspace{-0.25cm}
%%%%%%%%%%%%%%%%%%%%%%%%%%%%%%%%%
%
\subsection{Adjoint-based method}\label{ss:abm}
%
%%%%%%%%%%%%%%%%%%%%%%%%%%%%%%%%%

We first introduce a numerical discretization of the forward problem. Then we detail the adjoint method in Section~\ref{sec:ADJ:Problem} for computing the gradients required when solving \eqref{eq:objfunction-bayesian}.

%%%%%%%%%%%%%%%%
%
\subsubsection{The forward problem\label{sec:FWD:Problem}}
%%%%%%%%%%%%%%%%

We represent \eqref{eq:powersystem} compactly by 
\begin{align}
  \label{eq:powersystem:Compact}
  M \dot{\zz} &= F(t,\zz;\ipar)\,,~ \zz(0)=[x(0),y(0)]^T\,,
\end{align}
where $\vec u := (\vec x,\vec y)$ denotes the state variables,
$F=[h(\cdot),g(\cdot)]^T$, and $M$ is the DAE mass matrix, which is
block identity for $\vec x$ variables and zero in the rest.  Note
that $M$ should not be confused with the parametric inertia
$\ipar$. %\cM{This is confusing. Do you invert for anything other than
%  the M parameters in the equation on the left? Note that
%  you are assuming M constant. If you have already divided
%  by M to the right and those are the m's you should be clear about that, and call
%  M something else.}~\cE{Mihai, M is the mass matrix in the sense of
%  DAEs and is not related to the parameters m. Are they too close and
%  create confusion or am I missig something? I attempted to clarify
%  this in the text.}
Equation~\eqref{eq:powersystem:Compact} is discretized by using a time-stepping method. For instance, a trapezoidal-rule discretization leads to
%
%\begin{align}
%  \label{eq:powersystem:Compact:Discr}
%\nonumber
$  M \zz_{\tidx+1} = M \zz_{\tidx} + \frac{\dt}{2} \left(F(t_{\tidx},\zz_{\tidx};\ipar) +
  F(t_{\tidx+1},\zz_{\tidx+1};\ipar)\right)$, 
%\end{align}
%
where $\dt=t_{\tidx+1}-t_{\tidx}$.
With fixed $\zz(0)$, each choice of the parameters $\ipar$ will
generate a new trajectory.

\subsubsection{The adjoint problem and gradient computation\label{sec:ADJ:Problem}}

%Starting with an initial guess for the parameter $\ipar$, the
%optimization method iteratively updates this parameter using gradient
%information of the cost $\mathcal J$ in
%\eqref{eq:objfunction-bayesian} with respect to $\ipar$.
%Since
%evaluating $\mathcal J$ requires the solution of the forward model
%(\eqref{eq:powersystem}), which depends on $\ipar$,
%derivatives of $\mathcal J$ with respect to $\ipar$ need to take into account
%this implicit dependence as well. Here, we express the gradient using
%the so-called adjoint method.
To facilitate the gradient computation needed to
solve~\eqref{eq:objfunction-bayesian}, we use a Lagrangian approach
that augments
%the negative log posterior (i.e, the
%regularized data misfit functional),
$\mathcal{J}$ with additional
terms consisting of the forward DAE problem~\eqref{eq:powersystem}. 
%
%This strategy is used to compute the  gradient needed by the
%optimization problem,
%$\nabla_{\ipar} \J(\ipar)$.
Using the discrete adjoint approach \cite{Hager_2000,Hong_2015}
in PETSc~\cite{petsc-web-page}, we obtain the following
discrete adjoint equations:%~\cite{Hager_2000,Hong_2015}:
%Every time the optimizatin
%strategy requires the gradient value for a given $\ipar$, the system
%is integrated forward using \eqref{eq:powersystem:Compact:Discr}, the
%solution is stored or checkpionted and then following
%\cite{Hager_2000,Hong_2015}, the discrete adjoint equations take the
%following form:
%
\begin{align}
\nonumber
  M^T \lambda^* & = \lambda_{\tidx+1} + \frac{\dt}{2} \left(F_\zz^T(\zz_{\tidx+1})
  \lambda^* + r_\zz^T(t_{\tidx+1},\zz_{\tidx+1})\right) ,\\
  \label{eq:theta:adj}
  \lambda_\tidx & = M ^T \lambda^* + \frac{\dt}{2}  \left(F_\zz^T(\zz_\tidx)
  \lambda^* +  r_\zz^T(t_\tidx,\zz_\tidx)\right)\\
  \nonumber
  \mu_\tidx &= \mu_{\tidx+1} + \frac{\dt}{2} \left(
  F_{\ipar}^T(\zz_{\tidx+1}) +  F_{\ipar}^T(\zz_\tidx) \right) \lambda^* +  \\
  \nonumber
   & \qquad\frac{\dt}{2}  \left( r_\ipar^T(t_{\tidx+1},\zz_{\tidx+1}) +
  r_\ipar^T (t_\tidx,\zz_\tidx) \right), 
\end{align}
with $\tidx= N-1, \dots, 0$, $\lambda_{N} = \vec{0}$, $\mu_{N}
=\vec{0}$, where $r = -\log(\like(\obs | \ipar))$, and the gradients
are defined by $F_\zz=\frac{\partial F}{\partial \zz}$,
$F_{\ipar}=\frac{\partial F}{\partial \ipar}$, $r_\zz=\frac{\partial
  r}{\partial \zz}$, and $r_{\ipar}=\frac{\partial r}{\partial \ipar}$.
%\zC{What's the difference between $\lambda^*$ and $\lambda_k$?}
%Emil: Lambda star is a temporary variable.

The gradient of $\mathcal{J}$ with respect to $\ipar $ can be found by enforcing that the
derivative of the Lagrangian with respect to $\zz$ and the adjoint variables $(\lambda, \mu)$ vanish. This is given
by
%$\nabla_{\ipar} \J(\ipar)=\mu_0-\nabla_\ipar
%\frac{1}{2}\|\ipar-\mpr\|^2_{\matrix{\Gamma}_{\text{prior}}^{-1}}$.
$\nabla_{\ipar} \J(\ipar)=\mu_0-\prcov^{-1}(\ipar-\mpr)$.

The iterative procedures of a gradient computation are as follows. First, given a parameter sample $\ipar$, DAE~\eqref{eq:powersystem:Compact} is
solved for the forward solution $\zz$. The solution $\zz$ is
stored or checkpointed and further used to evaluate the data misfit
term $r$ in~\eqref{eq:theta:adj}. The adjoint equation is then solved (backward in time) to obtain the adjoint
solution $(\lambda, \mu)$.  Both the forward and adjoint solutions, along
with the current parameter  $\ipar$, are used to
evaluate $\nabla_{\ipar} \J(\ipar)$. Thus, a gradient computation requires two (forward and adjoint) DAE solves.
%
%The adjoint is evaluated along forward trajectory (need to compute the
%forward model). We need to do one forward and one adjoint computation
%for each gradient evaluation. 
%
%The inverse problem consists in finding the parameter values that
%support the observed data in the sense of \eqref{eq:posterior}. In
%this study we esimate the maximum a posteriori estimate (MAP) point of
%the posterior distribution \eqref{eq:posterior}. We compute the MAP
%point by using an adjoin-based optimization method.  To solve the
%forward and adjoint problems, required by the gradient evaluation, we
%use the time-stepping (TS) component of PETSc that provides ODE and
%DAE integrators. The optimization problem is then solved with the
%bounded limited memory variable metric quasi-Newton method for
%nonlinear minimization with bound constraints implemented in
%TAO~\cite{MunsonSarichWildEtAl12}.
%%a gradient-based optimization method (in TAO \cm{expand
%%  this}). This requires the gradient of $\pi_{\text{post}}(\ipar)$
%%with respect to $\ipar$.  
%
%a gradient-based optimization method (in TAO \cm{expand
%  this}). This requires the gradient of $\pi_{\text{post}}(\ipar)$
%with respect to $\ipar$.  

\subsubsection{Numerical solution to posterior minimization}\label{ss:optt}
%\zC{I suggest the following description: With the above-described gradient computation, we solve~\eqref{eq:objfunction-bayesian} using a quasi-Newton
%method with bound constraints implemented in TAO~\cite{MunsonSarichWildEtAl12}. This solver does not need evaluating second-order derivatives~\cite{Nocedal_book}; instead, the Hessian is approximated based on the previous evaluations of $\mathcal J(\ipar)$ and $\nabla_{\ipar} \J(\ipar)$. This approach has an asymptotic superlinear convergence  rate, and it updates the solutions iteratively by performing a Mor\'{e}-Thuente search~\cite{More_1994_LineSearch}  along the quasi-Newton direction until $\|\nabla_{\ipar} \J(\ipar)\|$ is small enough. }~\nnote{I like the paragraph below better, it's technical and sharp.}
The optimization problem~\eqref{eq:objfunction-bayesian} is solved with the bounded limited-memory variable-metric quasi-Newton
method for nonlinear minimization  implemented
in TAO~\cite{MunsonSarichWildEtAl12}. The method maintains a secant approximation to the Hessian from a limited number of previous evaluations of $\mathcal J(\ipar)$ and $\nabla_{\ipar} \J(\ipar)$ and uses this approximation to compute the quasi-Newton search direction. This approach achieves asymptotic superlinear convergence characteristic of Newton
method, but without evaluating second-order derivatives~\cite{Nocedal_book}. The numerical estimation procedure starts with an initial guess for $\ipar$ and iteratively updates this parameter by performing a Mor\'{e}-Thuente search~\cite{More_1994_LineSearch}  along the quasi-Newton direction. During this search a couple of evaluations of $\mathcal J$ may be needed in order to ensure sufficient decrease. The process stops when $\|\nabla_{\ipar} \J(\ipar)\|$ is small, which indicates that $\ipar$ is a local minimizer. 

\myvspace{-0.25cm}
%%%%%%%%%%%%%%%%%
\subsection{Stochastic spectral method}\label{ss:ssm}
%
%%%%%%%%%%%%%%%%%

For the second method, the DAE \eqref{eq:powersystem:Compact} is simulated at a small
number of samples to build a surrogate model.
Then, the obtained surrogate model (instead of the forward
solver) is used in the subsequent optimization to estimate
the parameters. This method is
particularly useful when the dimension of
the parameter space is small and the forward solver has a large state-space dimension, 
because it saves on the number of forward dynamic
simulations. Our example has $3$ parameters and $21$ state variables, so it
belongs to this category.  

%\cM{To Zheng and Emil: I replaced 3 with n, and d with n, here the
%  dimension of the parameter space. This is what Noemi uses
%  in the first paragraph. See if I got it wrong}
Given the prior density function of $\ipar$, a set of polynomial chaos basis functions $\left \{ \Psi _{\vec{\alpha}}(\ipar)\right \}_{|\vec{\alpha}|=0}^p$ are specified. Here $\vec{\alpha}\in \mathbb{N}^n$ is an index vector, $|\vec{\alpha}|$ denotes the $\ell _1$ norm, and positive integer $p$ is the highest order of the basis functions. These basis functions are orthornormal to each other, namely
%\begin{equation}
%\label{uni_gPC}
$\int\limits_{\mathbb{R}^n}  {   \Psi _{\vec{\alpha}}(\ipar)  \Psi
  _{\vec{\beta}}(\ipar)  \prior( {\ipar } )d\ipar
}=\delta_{\vec{\alpha},\vec{\beta}}$. 
%\end{equation}
Then, $\vec {f} (\vec {m})$ is approximated by a truncated generalized polynomial chaos expansion
%\begin{equation}
%\label{surrogate_i}
$\ff (\ipar) \approx \hat{\ff} (\ipar)=
\sum\limits_{|\vec{\alpha} | \leq p} {\vec{c}_{\vec{\alpha}} \Psi
  _{\vec{\alpha}}  (\ipar)}$ 
%\end{equation}
with $\vec{c}_{\vec{\alpha}}$ defined as
%$\begin{equation*}
%\label{c_project}
$\vec{c}_{\vec{\alpha}}=\int\limits_{\mathbb{R}^n}
{\Psi_{\vec{\alpha}} ( \ipar ) \vec{f} ( \ipar ) \prior(\ipar) d\ipar
}$. % \approx \sum_{i=1}^N {w_i \Psi_{\vec{\alpha}} ( \ipar _i) \vec{f} ( \ipar _i)}.
%$\end{equation*}
The total number of basis functions is $K=(p+n)!/(p!\,n!)$.
%Since there is a one-to-one correspondence between
%$\vec{\alpha}$ and integer $k\in [1,K]$, we can also denote
%$ \Psi _{\vec{\alpha}}(\ipar)$ by  $\Psi _{k}(\ipar)$.

In our implementations, $\vec{c}_{\vec{\alpha}}$ are computed
in two ways. The first choice is to employ projection-based stochastic collocation~\cite{col:2005,Ivo:2007,Nobile:2008}. Let $\left\{ \ipar_i, w_i\right\}_{i=1}^N$ be a set of quadrature points and weights corresponding to a numerical integration rule in the parameter space. Then we have 
%\begin{equation*}
%\label{c_project}
$\vec{c}_{\vec{\alpha}} \approx \sum_{i=1}^N {w_i \Psi_{\vec{\alpha}} ( \ipar _i) \vec{f} ( \ipar _i)}$.
%\end{equation*}
Popular methods for choosing the quadrature points include
tensor-product rules and sparse-grid methods~\cite{Gerstner:1998}. The
former needs $(p+1)^n$ samples to simulate the dynamic power systems,
whereas the latter needs fewer samples by using nested grid
samples. The second way is to use an interpolation method such as stochastic
  testing~\cite{zzhang:tcad2013}. Specifically, $K$ samples
  are selected, and the $\vec{c}_{\vec{\alpha}}$'s
  are obtained by solving a linear
  equation. In~\cite{zzhang:tcad2013}, a set of samples are
  generated by a quadrature rule (such as a tensor product
  Gauss-quadrature method); then $K$ dominant samples
  $\{ \ipar_j\}_{j=1}^K$ are subselected such that the
  matrix $\boldmath{V}$ (with its $j$th row being made of
  $\Psi_{\vec{\alpha}}(\ipar_j)$, $|\vec{\alpha}|\leq p$) is well conditioned.

% Here $(w_i, \ipar_i)$ is a pair of quadrature weight and sample for numerically evaluating the integral.

With a $p$th-order polynomial chaos expansion for $\vec{f} (\ipar)$, the negative log posterior now becomes a non-negative $2p$th-order polynomial function. We first write it as a combination of polynomial chaos basis function by stochastic collocation, then convert it to the summation of monomials:
%\begin{equation*}
$- \log \post (\ipar) \approx \sum _{ |\vec{\alpha}| =0}^{2p} {  q_{\vec{\alpha}}  \ipar^{  \vec{\alpha}  }   }, \; {\rm with} \; \ipar^{\vec{\alpha}}=m_1^{\alpha_1}m_2^{\alpha_2}\ldots m_n^{\alpha_n}$.
%\end{equation*}
 With this surrogate model, \eqref{eq:objfunction-bayesian} is simplified to
%\begin{equation}
%\label{eq:objfunction-surrogate}
%$\vec{\map} = \underset{\ipar}{\argmin}
%\; \hat{\J}(\ipar) :=  \sum _{ |\vec{\alpha}| =0}^{2p} {
%  q_{\vec{\alpha}}  \ipar^{  \vec{\alpha}  }   }$. 
$\vec{\map} = \underset{\ipar}{\argmin} \sum _{ |\vec{\alpha}| =0}^{2p} {
  q_{\vec{\alpha}}  \ipar^{  \vec{\alpha}  }   }$. 
%\end{equation}
This nonconvex optimization can be solved locally with gradient-based
methods as in Section~\ref{ss:abm} or globally with specialized
polynomial optimization solvers such as
GloptiPoly~\cite{Lasserre_2001,Lasserre_2002, Lasserre_2003}, when the
parameter dimension is low.

\section{Numerical Results}\label{s:numerical}
%In our approach the products are the MAP estimate of the parameters
%and the covariance matrix. The underlying hypothesis is that this
%estimate is a good representation of the real parameters that
%generated the data and that the variance is a good representation of
%the error between the MAP estimate and the real data.
In this section we evaluate how well does the MAP estimate approach the
parameters used to generate the data and how good of an error indicator is
the posterior variance.
%Here, we quantify and test these assumptions, discuss their
%limitations, and posit operational ranges and circumstances under
%which they can be used. The key direction is to compare the difference
%between the parameter values that generate the data and the MAP
%estimates in absolute terms and in relationship to the spread of the
%Bayesian posterior.  We also discuss the computational features and
%requirements of the two methods of computing the MAP estimate that we
%propose.

\begin{figure}
  \centering
  \begin{tabular}{c}
    \includegraphics[width=0.49\columnwidth]{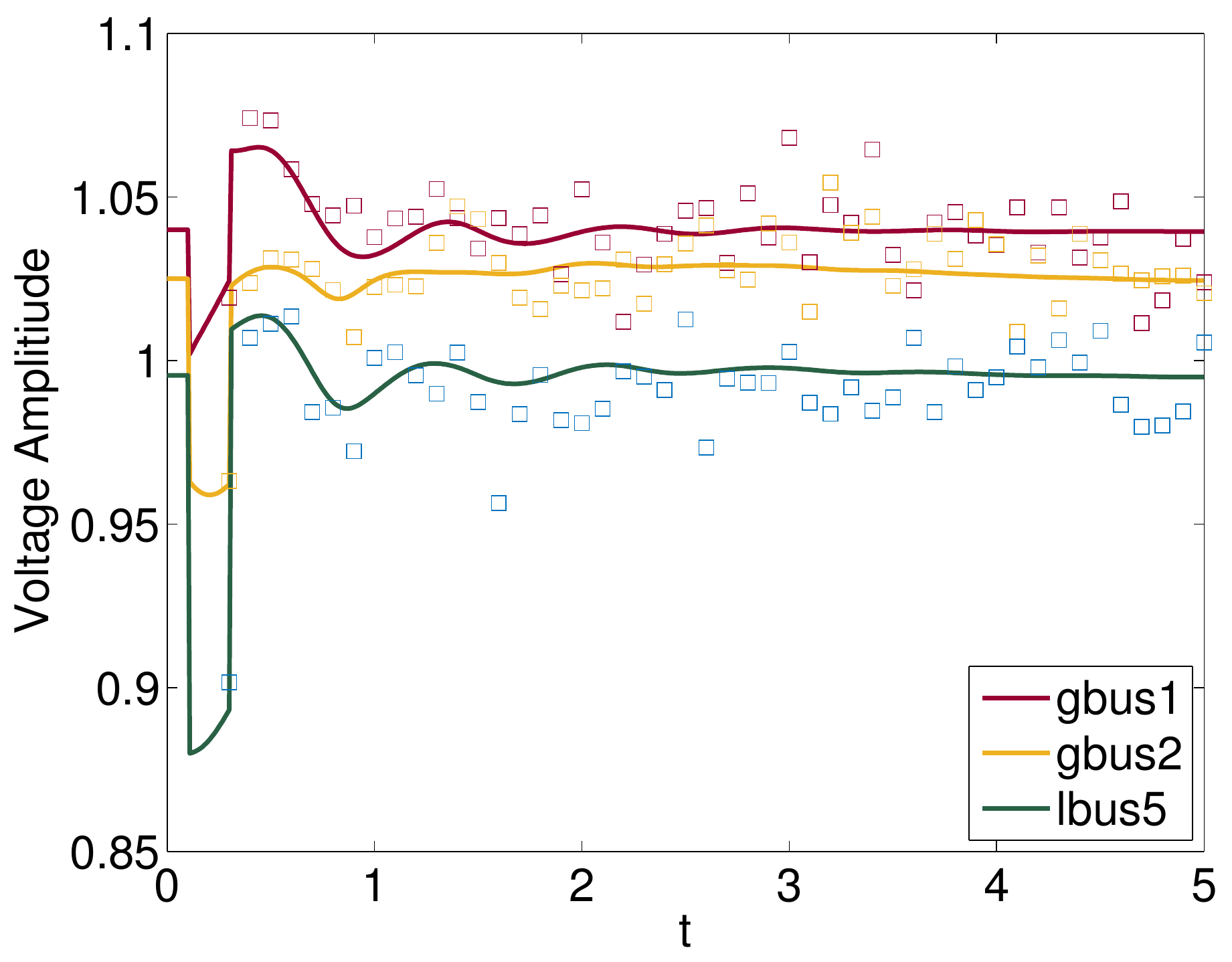}
    \includegraphics[width=0.47\columnwidth]{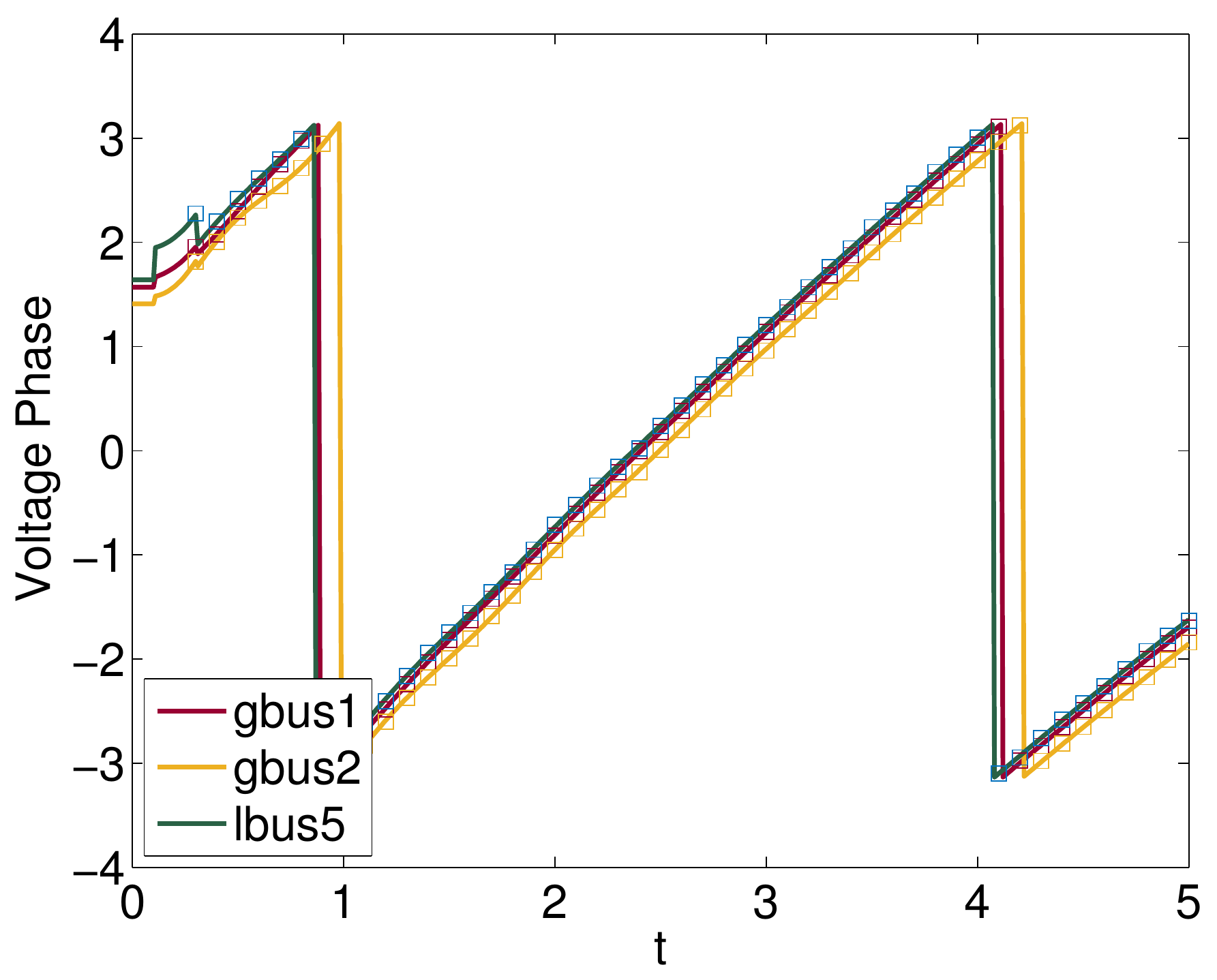}
  \end{tabular}
  \myvspace{-0.2in}
  \caption{Bus voltage amplitude (left) and phase (right) for the
    load buses 1 and 2 and generator bus 5 (red, yellow and green
    solid lines, respectively) based on the ``truth'' inertia
    parameter $\ipartrue = [23.64,6.40,3.01]$.
    %A load disturbance at $t
    %= 0.1$s was inserted to synthetically provoke a transient that
    %we measure. This was achieved by changing the load parameter
    %at the load bus 5 at $t = 0.1$s (from $1.25$ to $4.25$) and
    %then restore it at $t = 0.3$s.
    The squares show the
    corresponding synthetic observations with 1\% noise.}
  \label{fig:obs}
\myvspace{-0.25in}
\end{figure}

%The example we use is the IEEE 9-bus system
%~\nnote{We need to be %consistent with how we call the 3 gen, 9 bus system model.}
%depicted in Figure~\ref{fig:3gen9bus}.
A load disturbance at $t= 0.1$\,s, constant for $0.2$\,s, is inserted in the 9-bus model to provoke a
 transient. Its value during the switching action, $L$, is what
 characterizes this disturbance.
 We use
 as $\prcov$ a diagonal matrix
 with diagonal entries $[5.76, 0.36, 0.09]$. The prior mean and the
 ``truth'' inertia values are $\mpr = [24.00,6.00,3.10]$ and
 $\ipartrue = [23.64,6.4,3.01]$, respectively.  We carry out the forward
 simulation of the DAE \eqref{eq:powersystem} and we create synthetic
 voltage measurements at all 9 buses.  Here, we consider the case of
 independent observations; hence $\ncov$ is a diagonal matrix, with
 diagonal entries for all computations, unless otherwise specified,
 $10^{-4}$.
%the standard deviation is therefore 0.01%
%by adding
%independent identically distributed gaussian noise,consistent with our
%data model \eqref{equ:noise-model}.
The resulting voltage amplitudes,
phases, and synthetic measurements are depicted in Figure
\ref{fig:obs}.  This synthetic data is then used in the Bayesian
framework encapsulated in \eqref{eq:objfunction-bayesian}. 
We aim to quantify the estimation
error as a function of $L$ and the frequency of the
observations, which we assume consist of time series of the
voltages at all 9 buses, mimicking PMU data streams. 

%\cM{We need to be consistent when describing the example. I
%  suggest the introduction to be only about the abstraction
%  and mentioning the fact that it is the 9 bus system, but
%  that all details should appear in numerical section}

\myvspace{-0.25cm}
\subsection{Computational Setup}

The IEEE 9-bus example is implemented by using 
%the widely-used Portable,Extensible Toolkit for Scientific Computing (PETSc).
PETSc and is available as a part of the PETSc distribution. %\cM{For Emil: Link to download}
%~\nnote{Why do we need to specify the link in the text? Can we just refer to PETSc, for instance~\cite{petsc-web-page}?}
For
future, larger, examples, the setup has the advantage of having
intrinsic parallel capabilities \cite{petsc-web-page,petsc-user-ref}. %\cM{For Emil: Reference to PETSc}
The forward and adjoint problems needed by TAO for the numerical minimization of the posterior, as described in Section~\ref{ss:abm}, are set up and solved by using the PETSc
time-stepping library for DAEs. %They provide function and gradient evaluations to TAO to perform the numerical minimization of~\eqref{eq:objfunction-bayesian} as described in \S \ref{ss:optt}.

Both methods used in this paper will produce MAP estimates and their variance for
the inertia parameters $\ipar$ (a common notation for them in power
engineering literature is $H$~\cite{SauerPai98}).
%Since we also have
%the truth parameter values from which the data was simulated, we can
%now compare the deterministic numerical error, which relates the
%estimates to the real values but is unknown, with the a posteriori
%variances, which are computable.  When undertaking such comparisons,
%one can use various measures.
We compare the MAP estimate to the truth using a relative error metric:
\begin{align}
Err=\sqrt{\frac{1}{n} \smash \sum_{i=1}^n (\ipar(i)-\ipartrue(i))^2/\ipartrue(i)^2}.
\label{eq:err}
\end{align}
To facilitate interpretation we use a similar relative metric for the standard
deviation.
%(which, in classical statistical analysis, is
%typically discussed only in absolute terms, and not in
%relation to the mean).
Specifically, we use the following formula to normalize the square root of the trace
of the Hessian inverse ($\tau$):
\begin{align}
\tau=\sqrt{\smash\sum_{i=1}^n\postcov(i, i)/\ipartrue^2(i)}.\label{eq:std_var_norm}
\end{align}

Another statistic of interest is the positioning
of the real parameter in relation to the distribution. This
is not completely captured by the variance, because there could be
significant bias in the estimation. To this end, 
we compute the cumulative normal scores (CNS) $p$ for the
{\em actual  values} in relation to the Gaussian approximation of the
distribution:
\begin{align}
  p_i= \mbox{\bf erf} \Big[(\ipar(i)-\ipartrue(i)) / \sqrt{\postcov(i,i)}\Bigr].\label{eq:cdf}
\end{align}
Here ${\bf erf}$ is the {\em error function}, the cumulative
density of the standard normal. 
CNS are between $0$ and $1$ and indicate
how likely is that the real parameters are drawn from
the aposteriori distribution, with the distinction that
values very close to either $0$ or $1$ are considered {\em unlikely}. 

To determine whether our analysis had a good outcome, we use
the following considerations. If the estimation
procedure is successful, the error $Err$ should be small by
engineering standards (a few percentages or less). If the
stochastic model is a good depiction of reality, then $\tau$
should be mostly larger than $Err$ but comparable. This
reflects the fact that we are uncertain about the parameter
used in the estimation (as opposed to the deterministic
case); but when the data is informative, the standard
deviations should be comparable to the error (though exact
relational statements are difficult). A 
measure of successful representation of the uncertainty
analysis and validation of the statistical approach is
that the standard confidence values contain the real
parameter. That is, the CNS of the real parameters should be
away from $0$ and $1$ (for example, in the $[0.1,0.9]$
range) but not clustered at $0.5$, which would indicate an
excessively conservative variance. 

\myvspace{-0.25cm}
\subsection{Results} \label{sec:par:results}

\subsubsection{Dependence on experimental design parameters} \label{sss:dedp}

Our approach has two experimental design parameters:
%There are two experimental design parameters in our approach. 
the length of the time horizon over which
the estimation is carried out and the frequency
of the data. We now present the behavior of
$Err$, $\tau$, and the CNS values as a function of the
various choices of these parameters. 

%\cC{This section needs to be revisited in the light of the new data. I can do this.}

\onetwotable{Table~\ref{tbl:study1_adj_pc}}{Tables~\ref{tbl:study1}
  and~\ref{tbl:study2}} shows the estimation results for different estimation horizons $t_f$ and data frequencies $\Delta_t^{\scriptsize{obs}}$, shown by the first column in (a) and (b), respectively.
%
%Our results are listed in
%\onetwotable{Table~\ref{tbl:study1_adj_pc}}{Tables~\ref{tbl:study1}
%  and~\ref{tbl:study2}}, which shows the estimation results for
%different final times (top) and frequency of observations (bottom).
%The first column ($t_f$ in (a) and $\Delta_t^{\scriptsize{obs}}$ in
%(b)) shows the final time and the measurement frequency, respectively;
The second, third, and forth columns ($m_i$, $i$ = 1, 2, 3) indicate
the inverse solution, that is, the MAP point obtained with the
adjoint-based and the surrogate-based methods, separated by ``$/$'';
the fifth column (\#iter) indicates the number of iterations taken by
the adjoint-based method to converge. The sixth and seventh columns
($\tau$ and $Err$) show the standard deviation normalized by the
``truth'' inertia parameter (as given in~\eqref{eq:std_var_norm}) for
the two methods and the deterministic error computed with the
adjoint-based method), respectively. The last
three columns show the $p$-values computed with the adjoint-based
method by using~\eqref{eq:cdf}. For these simulations the forward problem
time step was $\Delta_t = 0.01$, the load parameter (at load bus 5)
was $5.5$, and the iterations were terminated when the norm of the
gradient fell below $10^{-6}$.
%In the first 3 columns, we display the 3 parameter values (denoted
%here as $H$), estimated by the two methods from \S \ref{s:solmet} in
%the form ``adjoint method estimate / stochastic spectral element
%estimate''. For the stochastic spectral method we have used degree 3
%polynomials. %\cM{Zheng and Emil, verify me on that}
%In column 4 we
%display the number of iterations of the adjoint based method, keeping
%in mind that for fixed values of the experimental design parameters
%the effort is linear in the number of iterations. In column 5 we
%display $\tau$ and the unscaled variance, whereas in column 6 we
%display the deterministic error, $Err$. Finally in columns 6-8 we
%display the the CNS values for each parameter. 

As we see in \onetwotable{Table~\ref{tbl:study1_adj_pc}}{Tables~\ref{tbl:study1}
  and~\ref{tbl:study2}}(b),
for data frequency
of 10 measurements per second or better, the deterministic
error $Err$ never gets above 2\%. In that range, the scaled standard
deviation $\tau$ is 5\% or better and 
%The ratio of $\tau/Err$ is most of the times less than 5, except for one case where it is slightly larger than 10. 
the ratio of $\tau/Err$ is always less than $4.5$.
The CNS values are
comfortably within $[0.1,0.9]$. We also note that the error,
$Err$, does not significantly improve with finer
measurements, and it oscillates with the decrease of the data frequency. On the other hand, the scaled standard
deviation does improve as the additional data reduces the
impact of the prior.
We conclude that in the range of 10
measurements per second or better, by the standards
indicated above, our statistical approach is
successful.
%That is,  it
%both creates estimates that are within 2\% of the real
%value and provides a statistical relative error estimate of the
%same order of magnitude as the deterministic error (and
%always less than 5\%).
Moreover, when providing confidence
intervals based on our Bayesian framework, the $[0.1,0.9]$
confidence interval {\em always} contains the real value. 
We also note that 10 measurements per second is comfortably
within the capabilities of typical PMU data streams of 30
measurements per second.

In \onetwotable{Table~\ref{tbl:study1_adj_pc}}{Tables~\ref{tbl:study1}
  and~\ref{tbl:study2}}(a), we list the effect of the
length of the estimation interval on the estimation. 
We do this at a data frequency of $20$ Hz ($\Delta_t^{obs}=0.05$), which is within the sampling rates ($30$ to $0.033$ Hz) supported by PMUs.
%While
%we do this at the  fine data spacing of 100 measurements per
%second, our previous discussion indicates that the results
%will likely not be very different for 10 measurements per
%second, insofar the deterministic error is
%concerned.
%\cM{To Cosmin: This will be a point of contention. We should
%  ahve done 0.1, as we indicated that to be the sweet
%  spot. If that is not too hard I would advocate it, if no
%   oh well}.
We  observe that for estimation horizons of $1$\,s or longer,
our statistical approach is also successful.
%, with one
%possible degradation (that we will describe later). 
In that range, both the deterministic error and the statistical errors are
less than $4 \%$,  and their ratio is never more than 4.
%in effect appears to converge to values very close to 1 as the time interval gets longer. 
Also, both error indicators
are decreasing with longer time horizons, whereas the
deterministic error was relatively insensitive to
measurement frequency. 

The CNS values  for the larger inertia, $m_1$,  fit comfortably within the $[0.1,0.9]$ confidence
interval. The smaller inertias, however, are contained only
in the $[0.01,0.99]$ range for the very long estimation
horizons. This interval, which for normal distributions
is about 3 standard deviations left and right of the mean, is
not abnormally wide by statistical analysis standards. But
it does suggest that smaller inertias
are harder to estimate accurately relative to larger ones,
which is not altogether surprising. However, when seen in the light
of the small relative standard deviation (about 1\%), 
those confidence intervals will be tight from an
operational perspective.

Therefore, when having the choice of more
frequent observations or longer estimation intervals, the
latter appears to be more beneficial to the quality of the
estimation once we are in range of $10$ measurements per
second or better.
%But an interval of $1$\,s or longer certainly produces
%satisfactory statistical outcomes. 

\onetwotable{}{
%---------- t a b l e  -------------------
\begin{table}[!t]
      \centering
     \caption{The effect of the time horizon on the
        ability to recover the inertia parameter for the power grid
        inverse problem. The first column ($t_f$) shows the final
        time; the second, third and forth columns ($\ipar_i$, i = 1,
        2, 3) indicate the inverse solution, i.e., the MAP point; and
        the fifth column (\#iter) indicates the number of iterations
        for the optimization solver to converge. This table shows that
        the minimum time horizon to recover the inversion parameter is
        $[0, 0.7]$. For these simulations the forward problem time
        step size was $\Delta_t = 0.01$, and the iterations were
        terminated when the norm of the gradient was decreased by a
        factor of $10^{6}$.\cE{Move some of the stuff to text
          (especially the last parts.} }
      \label{tbl:study1}
      \begin{tabular}{|p{0.4cm}|p{0.8cm}|p{0.8cm}|p{0.8cm}|p{0.5cm}|p{1.2cm}|p{1.cm}<{$}}
      \hline
      \multirow{2}{*}{$t_f$}   &   {${m_1}$}   &    {$m_2$}   & {$m_3$} & \#iter & var\\ %error\\
      \cline{2-6}
      & \multicolumn{5}{c|}{$\Delta_t$ = 0.01, $\Delta_t^{\scriptsize{obs}}$ = 0.01} \\
      \hline
      5     &  23.63     &  6.34    &  3.05   &   9   & 0.003 \\%5.42e-03 \\
      2     &  23.71     &  6.35    &  3.08   &   11  & 0.007 \\%8.24e-03 \\
      1     &  23.70     &  6.30    &  3.11   &   13  & 0.029 \\%1.23e-02 \\
      0.9   &  23.87     &  6.35    &  3.09   &   12  & 0.082 \\%9.79e-03 \\
      0.8   &  23.79     &  6.35    &  3.07   &   10  & 0.194 \\%7.44e-03 \\
      \rowcolor[gray]{0.8}
      0.7   &  23.81     &  6.36    &  3.06   &   13  & 0.308 \\%6.38e-03 \\
      0.6   &  23.82     &  6.01    &  3.16   &   11  & 0.768 \\%2.64e-02 \\
      0.5   &  22.83     &  6.05    &  3.12   &   8   &  2.455\\%2.47e-02 \\
      \hline
      \end{tabular}
\end{table}

\begin{table}[!t]
  \centering
  \begin{tabular}{|p{0.7cm}|p{0.5cm}|p{0.9cm}|p{0.9cm}|p{0.9cm}|p{0.5cm}|p{1.2cm}|}
    \hline \multirow{2}{*}{$\Delta_t^{\scriptsize{obs}}$} &
    \multirow{2}{*}{$N_{obs}$} & {${m_1}$} & {\bf $m_2$} & {\bf
      $m_3$} & {\#iter} & var \\%error\\
    \cline{3-7} &&
    \multicolumn{5}{c|}{$t_f$ = 0.7, $\Delta_t$ = 0.01} \\
    %, $N_t$=61
    \hline
    0.01 & 41 & 23.81     &  6.36    &  3.06   &   13  & 0.3084 \\%6.38e-03 \\
    0.02 & 21 & 23.61     &  6.39    &  3.03   &   12  & 0.5345 \\%2.31e-03\\
    0.04 & 11 & 22.89     &  6.34    &  2.91   &   9   & 0.8548 \\%1.56e-02 \\
    \rowcolor[gray]{0.8}
    0.1  & 5  & 23.92     &  6.25    &  2.94   &   10 &  1.4934 \\%1.17e-02 \\
    0.2  & 3  & 24.44     &  6.37    &  2.95   &   10 & 2.0566\\%1.32e-02\\
    0.4  & 2  & 24.40     &  6.40    &  2.86   &   9  & 2.3758\\%1.98e-02 \\
    \hline
  \end{tabular}
  \caption{A study of the ability to recover the inertia parameter for
    the power grid inverse problem. We fix the final time to $0.6$s
    and let the frequency of the data vary as listed in the first
    column ($\Delta_t^{obs}$); the second column ($N_{obs}$) shows the
    number of observation points; the fourth, fifth and sixth columns
    ($\ipar_i$, i = 1, 2, 3) indicate the inverse solution; and the
    last column (\#iter) indicates the number of iterations for the
    optimization solver to converge. This table shows that a frequency
    of one observation at every $0.1s$ produces good parameter
    reconstruction. For these simulations, the forward problem time
    step size was $\Delta_t = 0.01$s, and the iterations were
    terminated when the norm of the gradient was decreased by a factor
    of $10^{6}$. }
  \label{tbl:study2}
\end{table}
}
%------------------------------------
 \onetwotable{
\begin{table*}[!t]
  \centering
  \caption{A study of the effect of the time horizon and
        measurement frequency on the ability to recover the inertia
        parameter for the power grid inverse problem. %The first column
        %($t_f$ in (a) and $\Delta_t^{\scriptsize{obs}}$ in (b)) shows
        %the final time and the measurement frequency,
        %respectively; the second, third and forth columns ($\ipar_i$,
        %$i$ = 1, 2, 3) indicate the inverse solution, i.e., the MAP
        %point obtained with the adjoint-based and the surrogate based
        %methods, separated by ``$/$''; the fifth column (\#iter)
        %indicates the number of iterations taken by the adjoint-based
        %method to converge. The sixth and seventh columns (std var and
        %error) show the standard deviation normalized by the ``truth''
        %inertia parameter (as given in~\eqref{eq:std_var_norm}) for
        %the two methods and the error computed with the adjoint-based
        %method~\eqref{eq:err}), respectively. Finally, the last three
        %columns show the p-values computed with the adjoint-based
        %method using~\eqref{eq:cdf}.
        %For these simulations the
        %forward problem time step size was $\Delta_t = 0.01$, the load
        %was $5.5$, and the iterations were terminated when the norm of
        %the gradient was decreased by a factor of $10^{6}$.
      }
  \begin{tabular}{|p{0.4cm}|p{1.5cm}|p{1.5cm}|p{1.5cm}|p{0.5cm}|p{2.2cm}|p{1cm}|p{0.9cm}|p{0.9cm}|p{0.9cm}|}
      \hline
      \multirow{2}{*}{$t_f$} & {${m_1}$} & {$m_2$} & {$m_3$} & \#iter &  $\tau$  & $Err$ & $p_1$ & $p_2$ & $p_3$\\
      \cline{2-10}
      & \multicolumn{9}{c|}{(a) $\Delta_t$ = 0.01, $\Delta_t^{\scriptsize{obs}}$ = 0.05} \\
      \hline

       5.0 & 23.60 / 23.60 & 6.35 / 6.37 & 3.02 / 3.00 & 15 & 1.59e-02 / 1.79e-02 &  5.17e-03 & 0.1679 & 0.1484 & 0.5804 \\
%       4.0 & 23.76 / 23.77 & 6.31 / 6.31 & 3.13 / 3.13 & 16 & 1.74e-02 / 2.11e-02 &  2.38e-02 & 0.9870 & 0.0564 & 0.9949 \\
       3.0 & 23.79 / 23.81 & 6.39 / 6.41 & 3.06 / 3.05 & 9  & 1.85e-02 / 2.04e-02 &  1.01e-02 & 0.9659 & 0.4163 & 0.8470 \\
%       2.0 & 23.53 / 23.53 & 6.36 / 6.39 & 3.00 / 2.98 & 11 & 2.31e-02 / 2.36e-02 &  5.02e-03 & 0.2182 & 0.2939 & 0.4048 \\
       1.0 & 23.56 / 23.55 & 6.32 / 6.32 & 3.06 / 3.06 & 14 & 3.60e-02 / 3.66e-02 &  1.30e-02 & 0.4019 & 0.2384 & 0.7423 \\
%       0.9 & 23.38 / 23.38 & 6.29 / 6.31 & 3.11 / 3.09 & 13 & 4.64e-02 / 4.65e-02 &  2.34e-02 & 0.3027 & 0.2546 & 0.8647 \\
       0.8 & 23.67 / 23.63 & 6.54 / 6.53 & 2.95 / 2.95 & 11 & 5.81e-02 / 5.80e-02 &  1.76e-02 & 0.5123 & 0.7314 & 0.2295 \\
%       0.7 & 24.18 / 23.17 & 6.08 / 6.07 & 3.08 / 3.08 & 11 & 6.76e-02 / 6.69e-02 &  3.41e-02 & 0.3215 & 0.0958 & 0.7401 \\
       0.6 & 22.45 / 22.45 & 6.14 / 6.13 & 3.01 / 3.00 & 10 & 9.43e-02 / 9.29e-02 &  3.74e-02 & 0.1924 & 0.2337 & 0.4892 \\
%       0.5 & 25.19 / 24.18 & 6.01 / 6.01 & 3.12 / 3.12 & 8  & 1.29e-01 / 1.29e-01 &  4.31e-02 & 0.6085 & 0.2283 & 0.7492 \\
      %\hline
      \cline{1-10}
      $\Delta_t^{\scriptsize{obs}}$ & \multicolumn{9}{c|}{ (b) $t_f$ = 1, $\Delta_t$ = 0.01}\\
      \hline
      0.01 & 23.25 / 23.23 & 6.29 / 6.28 & 2.97 / 2.97 & 12 & 1.65e-02 / 1.67e-02 &  1.56e-02 & 0.0078 & 0.0195 & 0.1641 \\
      0.02 & 23.81 / 23.76 & 6.50 / 6.49 & 3.00 / 2.99 & 12 & 2.34e-02 / 2.36e-02 &  9.79e-03 & 0.7635 & 0.8897 & 0.4218 \\
      0.05 &  23.56 / 23.55 & 6.32 / 6.32 & 3.06 / 3.06 & 14 & 3.60e-02 / 3.66e-02 &  1.30e-02 & 0.4019 & 0.2384 & 0.7423 \\
      0.10 & 22.91 / 23.87 & 6.42 / 6.43 & 3.06 / 3.04 & 13 & 4.82e-02 / 4.89e-02 &  1.11e-02 & 0.7332 & 0.5607 & 0.6522 \\
      0.35 & 23.53 / 23.52 & 6.23 / 6.21 & 2.98 / 2.98 & 11 & 9.04e-02 / 9.08e-02 &  1.69e-02 & 0.4297 & 0.2452 & 0.4412 \\
      %\hline
      %0.01 & 23.81 / 23.82 &  6.36 / 6.36 &  3.06 / 3.06   &   13  & 3.74e-02 / 0.295 &  1.19e-02 & 0.6259 & 0.4021 & 0.8113 \\
      %0.02 & 23.61 / 23.65 &  6.39 / 6.40 &  3.03 / 3.03   &   12  & 4.98e-02 / 0.514 &  4.74e-03 & 0.4884 & 0.4884 & 0.6123 \\
      %0.04 & 22.89 / 22.91 &  6.34 / 6.34 &  2.91 / 2.91   &   9   & 6.34e-02 / 0.833 &  2.63e-02 & 0.2001 & 0.4124 & 0.1886 \\
      %\rowcolor[gray]{0.8}
      %0.1  & 23.92 / 23.94 &  6.25 / 6.26 &  2.94 / 2.94   &   10  & 8.30e-02 / 1.439 &  1.96e-02 & 0.5958 & 0.3118 & 0.3196 \\
      %0.2  & 24.44 / 24.44 &  6.37 / 6.37 &  2.95 / 2.96   &   10  & 9.69e-02 / 2.014 &  2.22e-02 & 0.7206 & 0.4721 & 0.3819 \\
      %0.4  & 24.40 / 24.39 &  6.40 / 6.40 &  2.86 / 2.87   &   9   & 1.27e-01 / 2.303 &  3.32e-02 & 0.6992 & 0.5076 & 0.3107 \\ 
      \hline
      \end{tabular}
      \myvspace{-0.2in}
      \label{tbl:study1_adj_pc}
\end{table*}
}{}

\onetwotable{}{
%---------- t a b l e 3  PC-based MAP-------------------
\begin{table}[!t]
      \centering
            \begin{tabular}{|p{0.9cm}|p{0.9cm}|p{0.9cm}|p{0.9cm}|p{0.9cm}|p{0.9cm}|}
      \hline
      \multirow{2}{*}{$t_f$}   &   {${m_1}$}   &    {$m_2$}   &  {$m_3$}  &var \\
      \cline{2-5}
      & \multicolumn{4}{c|}{$\Delta_t$ = 0.01, $\Delta_t^{\scriptsize{obs}}$ = 0.01} \\
      \hline
%      5     &  23.64     &  6.41    &  2.99   &   19     \\
%      %2     &  23.62     &  6.43    &  2.98   &   14     \\
%      %1.5   &  23.62     &  6.43    &  2.98   &   13     \\
%      1     &  23.62     &  6.43    &  2.97   &   14     \\
%      %0.9   &  23.63     &  6.44    &  2.97   &   10     \\
%      0.8   &  23.59     &  6.45    &  2.98   &   10     \\
%      %0.7   &  23.60     &  6.59    &  2.99   &   11    \\
%      \rowcolor[gray]{0.8}0.6   &  23.60     &  6.44    &  3.00   &   12    \\
%      0.5   &  23.87     &  6.36    &  3.04   &   12     \\
%      0.4   &  23.88     &  6.43    &  2.98   &   13     \\
      %      \hline
      5     &  23.66     &  6.36    &  3.05    &  0.003  \\
      2     &  23.73     &  6.36    &  3.08   &  0.007  \\
      1     &   23.70    &    6.31   &   3.11   & 0.029 \\
      0.9   &    23.86    &    6.36   &    3.09   &  0.078 \\
      0.8   &    23.78    &   6.36    &    3.07   &  0.186 \\
      \rowcolor[gray]{0.8}
      0.7   &    23.82    &   6.36    &      3.06  &  0.295 \\
      0.6   &   23.84     &   6.02    &      3.17  &  0.734\\
      0.5   &    22.86    &   6.06    &     3.13    & 2.393 \\
      \hline
      \end{tabular}
            \caption{Solution from the 2nd-order
              polynomial-chaos-based MAP algorithm with the setting of
              Tab. \ref{tbl:study1}.}
            \myvspace{-0.2in}
      \label{tbl:study1_PC}
\end{table}

\begin{table}[!t]
      \centering
     \begin{tabular}{|p{0.9cm}|p{0.9cm}|p{0.9cm}|p{0.9cm}|p{0.9cm}| p{0.9cm}|}
      \hline \multirow{2}{*}{$\Delta_t^{\scriptsize{obs}}$} &
      \multirow{2}{*}{$N_{obs}$} & {${m_1}$} & {\bf $m_2$} & {\bf
        $m_3$} & var \\
      \cline{3-6}
      && \multicolumn{4}{c|}{$t_f$ = 0.7, $\Delta_t$ = 0.01} \\
      %, $N_t$=61
      \hline
%      0.01 & 60 & 23.69 & 6.39 & 3.01 & 10 \\
%      0.05 & 30 & 23.60 & 6.44 & 3.00 & 12 \\
%      \rowcolor[gray]{0.8}0.1 & 6   & 23.59 & 6.41 & 3.04 & 13 \\
%      0.5 & 2   & 23.20 & 6.34 & 3.53 & 9 \\
      0.01 & 41 & 23.82     &  6.36    &  3.06 &  0.295  \\
      0.02 & 21 & 23.65     &  6.40    &  3.03  &  0.514\\
      0.04 & 11 &22.91      &  6.34    &   2.91  & 0.833 \\
      \rowcolor[gray]{0.8}
      0.1  & 5  &  23.94     &   6.26   &   2.94  &1.439 \\
      0.2  & 3  &   24.44    &   6.37    &  2.96  & 2.014\\
      0.4  & 2  &    24.39   &  6.40l    &   2.87  & 2.303\\
      \hline
    \end{tabular}
    \caption{Solution from the 2nd-order polynomial-chaos-based MAP algorithm with the setting of
              Tab. \ref{tbl:study2}. }
    \label{tbl:study2_PC}
    \myvspace{-0.2in}
\end{table}
}
%%%%%%%%%%%%%%%%%%%%%%%%%%%%%%%%%
%
% \subsection{Statistical analysis}
%
%%%%%%%%%%%%%%%%%%%%%%%%%%%%%%%%%

\subsubsection{Dependence on the nature of the perturbation}

%We now investigate the dependence of the statistical
%quantities on the nature of the perturbations, of the
%problem, which are the size $L$ of the load perturbation and
%the standard deviation of the measurement noise
%$\sigma_m$. 
Having established in Section~\ref{sss:dedp} that the Bayesian
posterior standard deviation is a good indicator of the
parameter error, we estimate its behavior with
the size of the load perturbation $L$. We note that if 
there were no perturbation, the system would be in
steady state, and its inertias would thus not be observable. We thus anticipate
that a larger perturbation would result in better estimation
properties and thus lower posterior variances.
%Because measurement noise is indicative of lack of information, we
%anticipate that larger $\sigma_m$ will result in larger
%posterior variances.  

\begin{figure*}
\centering
\begin{tikzpicture}
  \node (l) at (-11,0.0)
  {\includegraphics[width=.6\columnwidth]{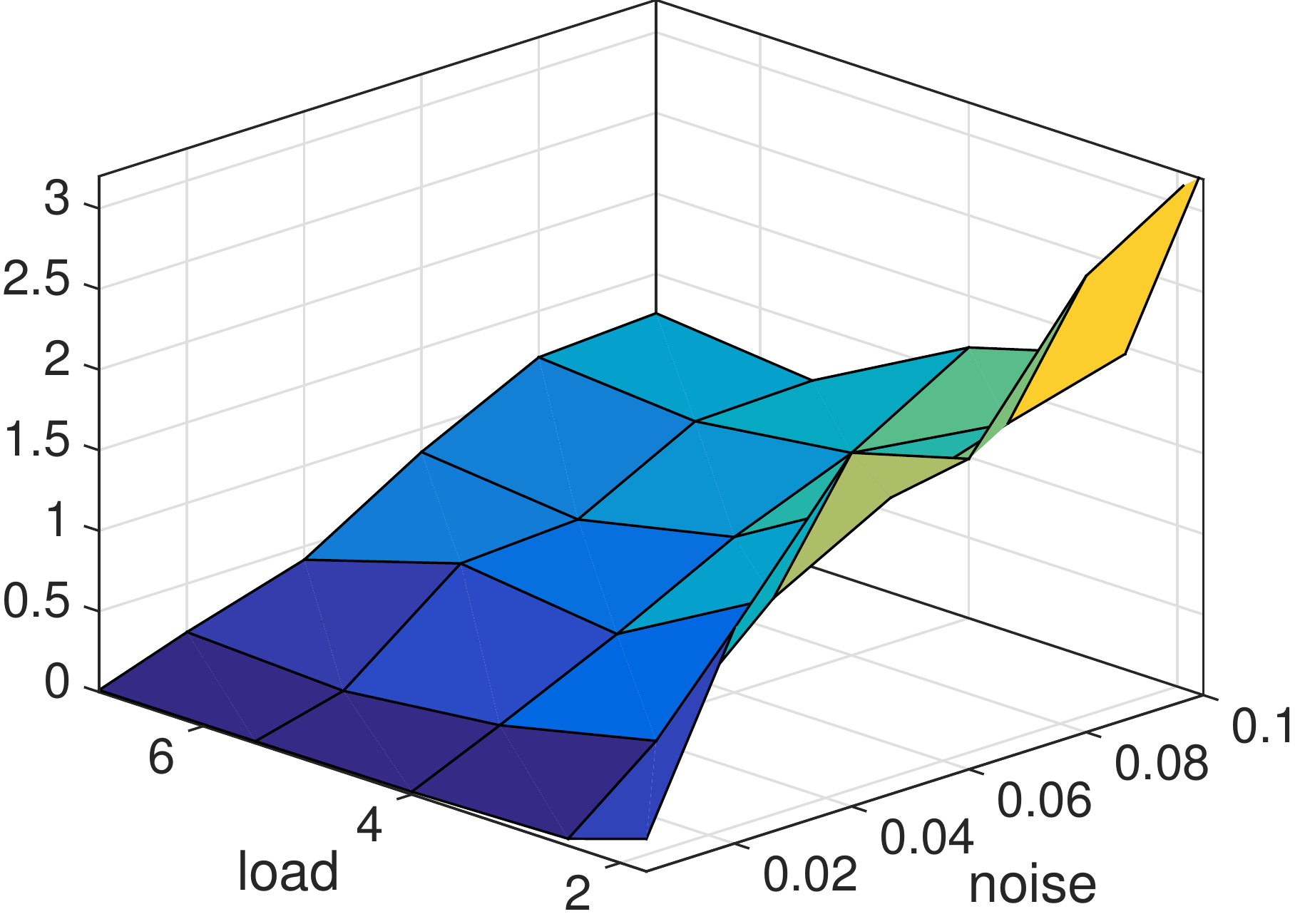}};
  %\node (l2) at (-5.5,0.5)
  %{\includegraphics[width=.6\columnwidth]{inertia_max_eig_tf2_nrsamples1_pn1e-1_3d.pdf}};
  \node (c) at (-8,0.1)
  {\includegraphics[width=.05\columnwidth]{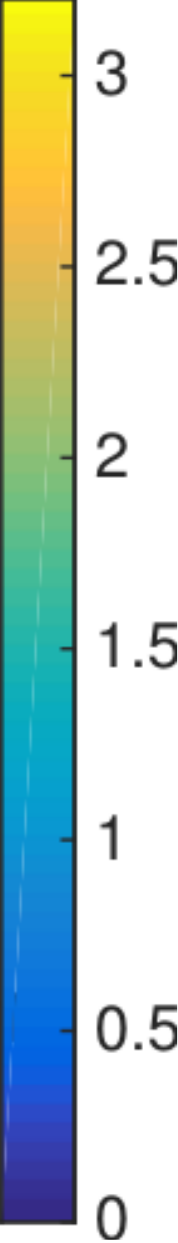}};
  \node (r) at (-2.0,0.0)
 {\includegraphics[width=1.2\columnwidth]{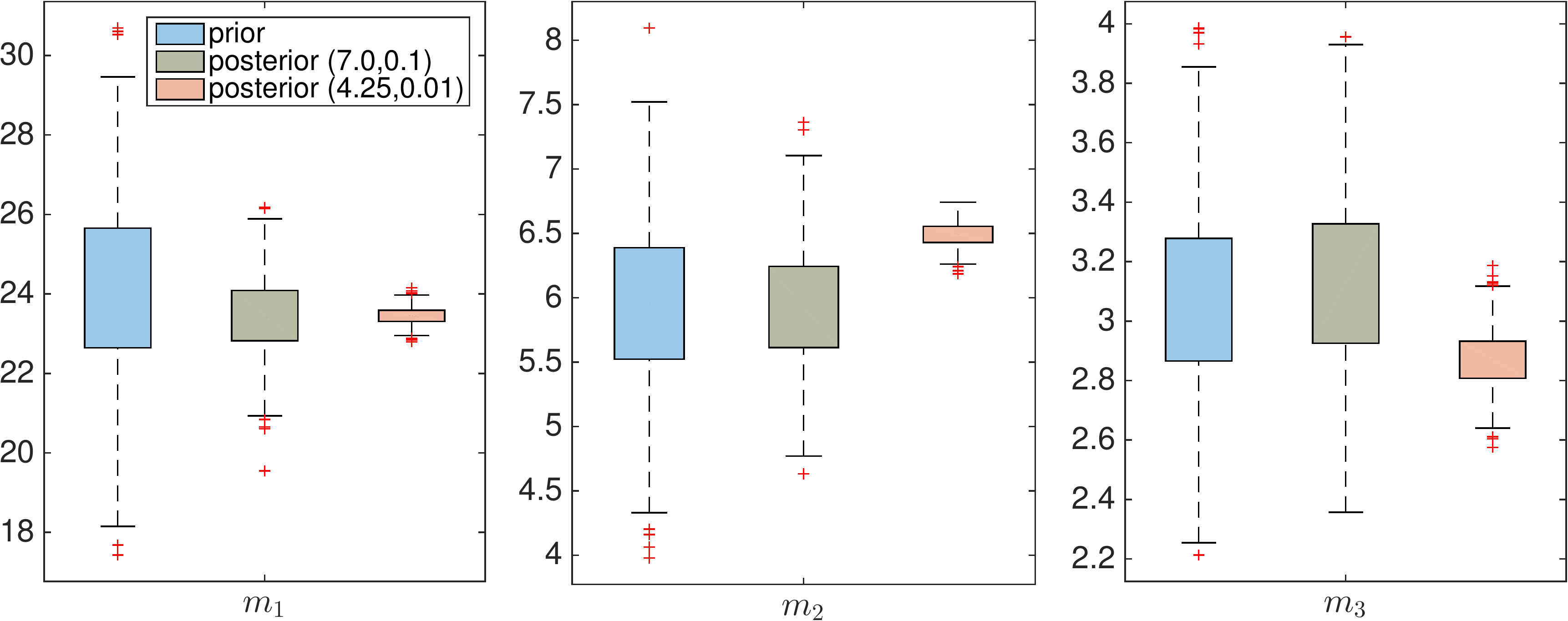}};
\end{tikzpicture}
\myvspace{-0.20in}
\caption{Left: Surface plot of the trace  of the Gaussianized
  posterior covariance as a function of noise and load for $t_f =
  2$ s.
  %This plot suggests that the variance increases as the noise
  %increases and the load decreases.
  Right: A ``whiskers boxplot'' of
  the prior and posterior mean and variances for load and noise values
  (7,0.1) and (4.25, 0.01) for the three inertia parameters.  The
  central mark is the median, the edges of the box are the $25$th and
  $75$th percentiles and the ``whiskers'' extend to the most extreme
  data points.}
\label{fig:tracecov}
\myvspace{-0.25in}
\end{figure*}

In Figure~\ref{fig:tracecov} we show a surface plot of the
trace of the Gaussianized posterior covariance (the sum of
the parameter variances) for several noise and
load values (left) and the ``whiskers boxplot'' of the prior and
posterior mean and variances for $L$  and $\sigma_m$ values
of (4.25, 0.01)
and (7,0.1), respectively.
%\cM{To Noemi and Cosmin: The whiskers boxplot is
%  confusing. It may have made sense to display all 4
%  choices, as it is nto clear if lower variance is due to
%  large load or large noise}
These figures show that, as anticipated, the variance increases as
the noise increases and the perturbation decreases, which indicates
that the deterministic error will have a similar behavior. 
The computational cost for computing MAP points (measured in number of
forward and adjoint solves) is shown in
Table~\ref{tbl:cost-adj-map}, which indicates that the
optimization effort is unaffected by the values of the
perturbation parameters.% $\sigma_m$ and $L$.

%---------- t a b l e  -------------------
\begin{table}[!t]
  \centering
   \caption{Computational cost for computing the MAP point
        measured in number of forward and adjoint solves.   %The
        %iterations were terminated when the norm of the gradient was
        %decreased by a factor of $10^{6}$. For all tests we used
        %$t_f = 2$s, $\Delta_t$ = 0.01, and $\Delta_t^{\scriptsize{obs}}$ =
        %0.1.
      }
      \begin{tabular}{|p{0.8cm}|p{0.8cm}|p{0.8cm}|p{2cm}|p{1cm}|}
      \hline
      Load & Noise & $\#$ Iter & $\#$ fwd/adj Solves & Time (s)\\
      \hline
      4.25 & 0.01  & 10 & 13 & 28\\
      4.25 & 0.1   & 9 & 13  & 31\\
      7.00 & 0.01  & 9 & 12  & 30\\
      7.00 & 0.1   & 10 & 14  & 34\\
      \hline
      \end{tabular}
     \myvspace{-0.14in}
     
      \label{tbl:cost-adj-map}
\end{table}

\myvspace{-0.25cm}
%%%%%%%%%%%%%%%%%%%%%%%%%%%%%%%%%
%
\subsection{Computational analysis}
%
%%%%%%%%%%%%%%%%%%%%%%%%%%%%%%%%%

We now discuss the computational cost for the two methods presented in
this study.  The adjoint-based method requires the value of the full
nonlinear model and its gradient for each iteration.  Additional
iterations may be required in the line-search procedure. The number of
forward and adjoint solves for selected cases is listed in
Table~\ref{tbl:cost-adj-map}. For these simulations we used $t_f =
2$ s, $\Delta_t$ = 0.01 s, and $\Delta_t^{\scriptsize{obs}}$ = 0.1 s. The
iterations for these simulations were terminated when the norm of the
gradient fell below $10^{-6}$.  To compare the adjoint-based method
cost with the stochastic spectral method, we need to account for the
cost of computing the adjoint, which is roughly the same as in the
forward run. In addition to the computational time, however, the
stochastic spectral element method has the advantage of working {\em
  without sensitivity information}. Given the considerable amount of
legacy software for which adjoints would be labor-intensive
to implement, this could confer it a practical
advantage. Moreover, the variance can
be naturally estimated with no additional cost, whereas the adjoint-based approach would need either finite differences or second-order adjoints
to compute the covariance.
%\cM{We should probably emphasize the fact that the spectral
%  element does not need adjoints. This is a nontrivial
%  selling point for this audience} 

%When analyzing the performance of the stochastic spectral
%element method, we note from \S \ref{ss:ssm}  that the total number of basis function used in polynomial-chaos expansion
%is $K=(d+n)!/(n!d!)$. 
%Of the 3 methods we used to construct these approximations,
%we note that stochastic testing
%requires $K$ samples to obtain the polynomial-chaos expansion;
%sparse-grid stochastic collocation requires about $2^p K$ samples as
%$d$ increases; tensor-product stochastic collocation technique is the
%most expensive, requiring at least $(p+1)^d$ samples in total.  The
%number of model evaluations in one of the cases under consideration
%are listed in Tab. \ref{tbl:study4_PC}.

%Insofar the results of the stochastic spectral element,
In Table~\ref{tbl:study3_PC} we show the costs of constructing
surrogate models using different approaches. In
Table~\ref{tbl:study4_PC} we show the MAP results using
different orders of surrogate models constructed by different
methods.  \onetwotable{}{For the 9-bus system, we use the configuration of Row 1 in
Table~\ref{tbl:study1} as a demonstration.}
Clearly, the
accuracy is significantly improved when we increase the order of
polynomial chaos expansion from 1 to 2, but the improvement is
marginal when we use third-order polynomial chaos
expansions. From these tables we see that we can obtain good-quality estimates of the parameters and their variance using
only 10 forward runs (degree 2). % when using the stochastic
%testing approach.
This is less intensive than the
adjoint-based method by a factor of about 2. 

\subsubsection{Challenges as we increase the number of parameters and
  the complexity of the problem}
The adjoint-based method has two major requirements: (1) code differentiation,
that is, the computation and implementation of derivatives such as the
ones in  \eqref{eq:theta:adj}, and (2) storing the forward trajectory
through checkpoints. 
Because only a few thousands of states need to be stored if the
time scales remain the same, even for interconnect size
examples this is unlikely to become a problem even on a
desktop. On the one hand point (1) is a significant
undertaking, although HPC tools such as PETSc increasingly
provide support for it natively. 
On the other hand, the
cost of the adjoint-based method is independent of the number of
parameters, and parallel implementations are also
possible.

The stochastic spectral method proved to be robust in
our experimental setting, requiring few model
evaluations to construct a viable surrogate. In addition,
all calculations can be trivially parallelized, and a
variance estimator is intrinsic. Moreover, once
a surrogate is obtained, one need not to
regenerate it if the model and setting do not
change. As discussed above, however, when the parameter
dimensionality $n$ is large, the number of simulation
samples required can be very large, leading to an extremely
high computational cost. Arguably one can obtain
efficiently a high-dimensional surrogate model by using some
advanced techniques, such as compressed
sensing~\cite{yxiu:2013}, tensor
recovery~\cite{zzhang:trecovery}, and proper generalized
decomposition~\cite{Nouy:2010}; but these techniques may still be inefficient for extremely high-dimensional cases (e.g., when $n>1000$). 

While a definitive comparison between the two approaches is difficult to make in
general because of the multiple features of the target problems, for a small number
of parameters and lack of sensitivity information, the
stochastic spectral element approach would be a strong
candidate for a solution. In our case, it did produce good
estimates a factor of 2 faster than the adjoint-based approach for
the proper choice of degree and construction method, which may be difficult to guarantee a priori. 

\subsubsection{Considerations about deployment}
While this is only an initial study, a practical implementation is worth considering. 
In such cases the initial state and the load would need to be
inverted as well. Because these are classical analyses, a
tiered approach is possible, where they are estimated
separately. One can, of course, create a
unified estimation approach with hybrid data sources; 
a mix of PMU and other data, such as SCADA, may need to
be considered.  While the performance of the method would
need itself to be re-evaluated, this can be done
in the Bayesian framework described in Section~\ref{s:probFor}.
As described, the method assumes that we have a way to
identify ``micro-transients'' suitable to trigger dynamical
estimation. This can be done for PMU
data. The method can also be modified to
support any type of perturbation, as well as in a
``rolling horizon'' approach, where it is not triggered but 
used continuously. This can be done, for
example, by restarting the estimation with the prior covariance
being the posterior one from the previous estimation
interval.
We anticipate that as long as 
the perturbations show 
enough dynamic range so that the method can excite transients that are
informative about the inertias, similar behaviors and
performance can be
expected. A more significant concern is the ability to
compute the estimate in real time. We note that forward
simulations for power grid transients using PETSc on
interconnect-sized networks have been run
faster than real time with less than 16 cores~\cite{Abhyankar_P2013}. %\cM{Emil, details
%from the work with Shri}.
  Therefore, for a few dynamic
parameters to invert with uncertainty, the stochastic spectral element method
could in principle work ``out of the box''.  For a large
number of dynamic parameters to invert, the issue is whether
the optimization can be fast enough. 
%While this cannot be answered positively at this time, 
Certainly a
promising direction is the usage of a rolling horizon
approach in conjunction with inexact optimization. 

\begin{table}[!t]
      \centering
  \caption{Total number of forward simulations to construct the surrogate models. }    
     \begin{tabular}{|c |c |c|c|}
      \hline \multirow{2}{*}{polyn. order} &
      \multicolumn{3}{c|}{Total number of forward simulations.} \\
            \cline{2-4}
      & stoch. testing & tensor prod. & sparse grid \\
      %, $N_t$=61
      \hline      1 & 4 & 8     &  7       \\
      2 & 10 & 27     &  19       \\
      3 & 20 & 64      &  39       \\
      \hline
    \end{tabular}
     \label{tbl:study3_PC}
     \myvspace{-0.2in}
\end{table}

\begin{table}[!t]
      \centering
  \caption{MAP results using different order of gPC expansions. }    
     \begin{tabular}{|c  |c | c|c|c|}
      \hline 
      surrogate model & gPC order & $m_1$ & $m_2$ & $m_3$ \\ \hline
        \multirow{3}{*}{stoch. testing}    
        & 1 &22.818  &6.745 & 2.248\\ \cline{2-5}
       & 2 &23.600  & 6.372&  3.000\\ \cline{2-5}
        & 3 &23.611  &6.351 &  3.021\\ \hline  
      \multirow{3}{*}{SC w/ tensor prod.}    
        & 1 & 23.751 & 6.420& 2.973\\ \cline{2-5}
       & 2 &  23.585& 6.374&  2.991\\ \cline{2-5}
        & 3 & 23.618 &6.347 & 3.026 \\ \hline    
        \multirow{3}{*}{SC w/ sparse grid}    
        & 1 &23.962  & 6.322& 3.156\\ \cline{2-5}
       & 2 & 23.584 & 6.375& 2.990 \\ \cline{2-5}
        & 3 & 23.617 &6.361 &  3.016\\ \hline
     \end{tabular}
     \myvspace{-0.2in}
    \label{tbl:study4_PC}
\end{table}

\section{Conclusions}

We have presented a Bayesian framework for parameter
estimation with uncertainty focused on the estimation of
dynamic parameters of energy systems. This investigation is
prompted by the rapid expansion of PMU sensors and the increased usage of renewable generation
whose inertia features may change in time and may not be
known to the stakeholder that must ensure transient stability
operation of the system. For such systems, inertia cannot be
assumed known and must thus be estimated together with its uncertainty. Because inertia has no
impact on steady-state features of the system, it needs transient scenarios under
which to be estimated. 

We have proposed two methods to compute the MAP estimates
and their variances: an adjoint-based method and a
stochastic spectral method. The former has the
advantage that it can compute gradients of the
log-likelihood function in a time that is a constant factor
of the one of the forward simulation irrespective of the
number of parameters considered. This method was implemented
in PETSc. The latter has the benefit
of needing no sensitivity capabilities, of employing only
forward simulations, and of providing an intrinsic estimate
of the variance.  It is suitable for the case of a
limited  number
of parameters. We have
demonstrated these methods on a 9-bus example case that is
available for download \cite{git_state_est}.
The three parameters to be estimated were the 
generator inertias. For this example we have generated synthetic
data of transient behavior by perturbing the load and adding
measurement noise
that we
have used to assess the behavior of our approaches. When
applying our method we have found that estimation time
horizons of $1$\,s or more and data frequency of at least 10
samples per second were sufficient for the error to be less
than 2\%, the posterior variance to be a good estimate of
the error, and some of the standard confidence intervals to
cover the real parameter (with the 3 standard deviations
ones always containing the real parameters). We have also observed that, as
expected, the error and posterior variance decrease with
increased system perturbation and decreased measurement
error. The computational effort was on the order of 10
forward simulations for the stochastic spectral method and 30 forward simulations for the adjoint-based method.
For usage in larger systems under real time constraints, and
under realistic data streams and use cases, further work may be necessary. 
Nevertheless, for the small parameter case the state of technology
is such that, with the use of parallel computing, the stochastic spectral method may
already provide sufficient capabilities.

\iffalse
\section*{Acknowledgements}
This material  was based upon work supported in part by the 
Office of Science, U.S. Dept. of Energy, Office of Advanced Scientific
Computing Research, under Contract DE-AC02-06CH11357. N.P. also acknowledges partial funding
through the U.S. Dept. of Energy, Office of Workforce Development for
Teachers and Scientists, Visiting Faculty Program.
\fi

\bibliographystyle{IEEEtran}
\bibliography{paperbib}

% Generated by IEEEtran.bst, version: 1.13 (2008/09/30)
\begin{thebibliography}{10}
\providecommand{\url}[1]{#1}
\csname url@samestyle\endcsname
\providecommand{\newblock}{\relax}
\providecommand{\bibinfo}[2]{#2}
\providecommand{\BIBentrySTDinterwordspacing}{\spaceskip=0pt\relax}
\providecommand{\BIBentryALTinterwordstretchfactor}{4}
\providecommand{\BIBentryALTinterwordspacing}{\spaceskip=\fontdimen2\font plus
\BIBentryALTinterwordstretchfactor\fontdimen3\font minus
  \fontdimen4\font\relax}
\providecommand{\BIBforeignlanguage}[2]{{%
\expandafter\ifx\csname l@#1\endcsname\relax
\typeout{** WARNING: IEEEtran.bst: No hyphenation pattern has been}%
\typeout{** loaded for the language `#1'. Using the pattern for}%
\typeout{** the default language instead.}%
\else
\language=\csname l@#1\endcsname
\fi
#2}}
\providecommand{\BIBdecl}{\relax}
\BIBdecl

\bibitem{ZhangBankWanEtAl13}
Y.~Zhang, J.~Bank, Y.-H. Wan, E.~Muljadi, and D.~Corbus, ``Synchrophasor
  measurement-based wind plant inertia estimation,'' in \emph{Green
  Technologies Conference, 2013 IEEE}.\hskip 1em plus 0.5em minus 0.4em\relax
  IEEE, 2013, pp. 494--499.

\bibitem{knyazkin2004parameter}
V.~Knyazkin, C.~Canizares, L.~H. S{\"o}der \emph{et~al.}, ``On the parameter
  estimation and modeling of aggregate power system loads,'' \emph{IEEE
  Transactions on Power Systems}, vol.~19, no.~2, pp. 1023--1031, 2004.

\bibitem{choi2006measurement}
B.-K. Choi, H.-D. Chiang, Y.~Li, H.~Li, Y.-T. Chen, D.-H. Huang, and M.~G.
  Lauby, ``Measurement-based dynamic load models: derivation, comparison, and
  validation,'' \emph{IEEE Transactions on Power Systems}, vol.~21, no.~3, pp.
  1276--1283, 2006.

\bibitem{bai2009novel}
H.~Bai, P.~Zhang, and V.~Ajjarapu, ``A novel parameter identification approach
  via hybrid learning for aggregate load modeling,'' \emph{Power Systems, IEEE
  Transactions on}, vol.~24, no.~3, pp. 1145--1154, 2009.

\bibitem{hiskens1999power}
I.~A. Hiskens and A.~Koeman, ``Power system parameter estimation,''
  \emph{Journal of Electrical and Electronics Engineering Australia}, vol.~19,
  pp. 1--8, 1999.

\bibitem{hiskens2001inverse}
I.~Hiskens, ``Inverse problems in power systems,'' \emph{Bulk Power System
  Dynamics and Control V}, pp. 180--195, 2001.

\bibitem{hiskens2004power}
I.~Hiskens \emph{et~al.}, ``Power system modeling for inverse problems,''
  \emph{IEEE Transactions on Circuits and Systems I: Regular Papers}, vol.~51,
  no.~3, pp. 539--551, 2004.

\bibitem{Tarantola05}
A.~Tarantola, \emph{Inverse Problem Theory and Methods for Model Parameter
  Estimation}.\hskip 1em plus 0.5em minus 0.4em\relax Philadelphia, PA: SIAM,
  2005.

\bibitem{SauerPai98}
P.~W. Sauer and M.~A. Pai, \emph{Power system dynamics and stability}.\hskip
  1em plus 0.5em minus 0.4em\relax Prentice Hall, 1998.

\bibitem{OakleyO'hagan07}
J.~E. Oakley and A.~O'hagan, ``Uncertainty in prior elicitations: a
  non--parametric approach,'' \emph{Biometrika}, vol.~94, pp. 427--441, 2007.

\bibitem{KaipioSomersalo05}
J.~Kaipio and E.~Somersalo, \emph{Statistical and Computational Inverse
  Problems}, ser. Applied Mathematical Sciences.\hskip 1em plus 0.5em minus
  0.4em\relax New York: Springer-Verlag, 2005, vol. 160.

\bibitem{Hager_2000}
W.~Hager, ``{R}unge-{K}utta methods in optimal control and the transformed
  adjoint system,'' \emph{Numerische Mathematik}, vol.~87, no.~2, pp. 247--282,
  2000.

\bibitem{Hong_2015}
H.~Zhang, S.~Abhyankar, E.~Constantinescu, and M.~Anitescu, ``Discrete
  sensitivity analysis of power system dynamics,'' \emph{draft}, 2016.

\bibitem{petsc-web-page}
\BIBentryALTinterwordspacing
S.~Balay, S.~Abhyankar, M.~F. Adams, J.~Brown, P.~Brune, K.~Buschelman,
  L.~Dalcin, V.~Eijkhout, W.~D. Gropp, D.~Kaushik, M.~G. Knepley, L.~C.
  McInnes, K.~Rupp, B.~F. Smith, S.~Zampini, and H.~Zhang, ``{PETS}c {W}eb
  page,'' 2015. [Online]. Available: \url{http://www.mcs.anl.gov/petsc}
\BIBentrySTDinterwordspacing

\bibitem{MunsonSarichWildEtAl12}
T.~Munson, J.~Sarich, S.~Wild, S.~Benson, and L.~C. McInnes, ``{TAO} 2.0 users
  manual,'' Mathematics and Computer Science Division, Argonne National
  Laboratory, Tech. Rep. ANL/MCS-TM-322, 2012.

\bibitem{Nocedal_book}
J.~Nocedal and S.~J. Wright, \emph{Numerical Optimization}, 2nd~ed.\hskip 1em
  plus 0.5em minus 0.4em\relax New York: Springer, 2006.

\bibitem{More_1994_LineSearch}
J.~J. Mor{\'e} and D.~J. Thuente, ``Line search algorithms with guaranteed
  sufficient decrease,'' \emph{ACM Trans. Math. Softw.}, vol.~20, no.~3, pp.
  286--307, Sep. 1994.

\bibitem{col:2005}
D.~Xiu and J.~S. Hesthaven, ``High-order collocation methods for differential
  equations with random inputs,'' \emph{SIAM J. Sci. Comp.}, vol.~27, no.~3,
  pp. 1118--1139, Mar 2005.

\bibitem{Ivo:2007}
I.~Babu\v{s}ka, F.~Nobile, and R.~Tempone, ``A stochastic collocation method
  for elliptic partial differential equations with random input data,''
  \emph{SIAM J. Numer. Anal.}, vol.~45, no.~3, pp. 1005--1034, Mar 2007.

\bibitem{Nobile:2008}
F.~Nobile, R.~Tempone, and C.~G. Webster, ``A sparse grid stochastic
  collocation method for partial differential equations with random input
  data,'' \emph{SIAM J. Numer. Anal.}, vol.~46, no.~5, pp. 2309--2345, May
  2008.

\bibitem{Gerstner:1998}
T.~Gerstner and M.~Griebel, ``Numerical integration using sparse grids,''
  \emph{Numer. Algor.}, vol.~18, pp. 209--232, Mar. 1998.

\bibitem{zzhang:tcad2013}
Z.~Zhang, T.~A. El-Moselhy, I.~M. Elfadel, and L.~Daniel, ``Stochastic testing
  method for transistor-level uncertainty quantification based on generalized
  polynomial chaos,'' \emph{IEEE Trans. CAD of Integr. Circuits and Syst.},
  vol.~32, no.~10, pp. 1533--1545, Oct 2013.

\bibitem{Lasserre_2001}
J.~B. Lasserre, ``Global optimization with polynomials and the problem of
  moments,'' \emph{SIAM J. Optimization}, vol.~11, no.~3, p. 796–817, 2001.

\bibitem{Lasserre_2002}
------, ``An explicit equivalent positive semidefinite program for 0-1
  nonlinear programs,'' \emph{SIAM J. Optimization}, vol.~12, no.~3, pp.
  756--769, 2001.

\bibitem{Lasserre_2003}
D.~Henries and J.~B. Lasserre, ``Gloptipoly: Global optimization over
  polynomials with {Matlab} and {SeDuMi},'' \emph{ACM Trans. Mathematical
  Software}, vol.~29, no.~2, p. 165–194, 2003.

\bibitem{petsc-user-ref}
\BIBentryALTinterwordspacing
S.~Balay, S.~Abhyankar, M.~F. Adams, J.~Brown, P.~Brune, K.~Buschelman,
  L.~Dalcin, V.~Eijkhout, W.~D. Gropp, D.~Kaushik, M.~G. Knepley, L.~C.
  McInnes, K.~Rupp, B.~F. Smith, S.~Zampini, and H.~Zhang, ``{PETS}c users
  manual,'' Argonne National Laboratory, Tech. Rep. ANL-95/11 - Revision 3.6,
  2015. [Online]. Available: \url{http://www.mcs.anl.gov/petsc}
\BIBentrySTDinterwordspacing

\bibitem{yxiu:2013}
X.~Yang and G.~E. Karniadakis, ``Reweighted $l_1$ minimization method for
  stochastic elliptic differential equations,'' \emph{J. Comp. Phys.}, vol.
  248, no.~1, pp. 87--108, Sept. 2013.

\bibitem{zzhang:trecovery}
Z.~Zhang, H.~D. Nguyen, K.~Turitsyn, and L.~Daniel, ``Probabilistic power flow
  computation via low-rank and sparse tensor recovery,'' \emph{arXiv preprint:
  arXiv:1508.02489}, pp. 1--8, Aug 2015.

\bibitem{Nouy:2010}
A.~Nouy, ``Proper generalized decomposition and separated representations for
  the numerical solution of high dimensional stochastic problems,'' \emph{Arch.
  Comp. Meth. Eng.}, vol.~27, no.~4, pp. 403--434, Dec 2010.

\bibitem{Abhyankar_P2013}
S.~Abhyankar, B.~Smith, and E.~Constantinescu, ``Evaluation of overlapping
  restricted additive {S}chwarz preconditioning for parallel solution of very
  large power flow problems,'' in \emph{Proceedings of the 3rd International
  Workshop on High Performance Computing, Networking and Analytics for the
  Power Grid}, ser. HiPCNA-PG '13.\hskip 1em plus 0.5em minus 0.4em\relax ACM,
  2013, pp. 5:1--5:8.

\bibitem{git_state_est}
\BIBentryALTinterwordspacing
 [Online]. Available:
  \url{https://github.com/Argonne-National-Laboratory/PowerSystemsEstimation.git}
\BIBentrySTDinterwordspacing

\end{thebibliography}
\myvspace{-0.35in}
\begin{IEEEbiography}[{\includegraphics[width=1in,clip,keepaspectratio]{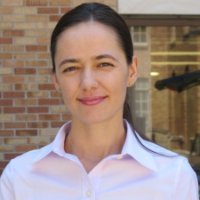}}]{Noemi Petra}
 is an assistant professor in the Applied Mathematics department in
 the School of Natural Sciences at the University of California,
 Merced. She earned her Ph.D. degree in applied mathematics from the
 University of Maryland, Baltimore County in 2010. Prior to joining UC
 Merced, Noemi was the recipient of the ICES Postdoctoral Fellowship
 in the Institute for Computational Engineering \& Sciences at the
 University of Texas at Austin. Her research interests include inverse
 problems, uncertainty quantification and optimal experimental design.
\end{IEEEbiography}
\myvspace{-0.3in}
\begin{IEEEbiography}[{\includegraphics[width=1in,clip,keepaspectratio]{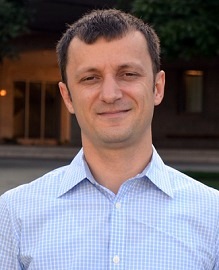}}\myvspace{-0.3in}]{Cosmin G. Petra} is an assistant computational mathematician in the Mathematics and Computer Science Division at Argonne National Laboratory. He received the B.Sc. degree in mathematics and computer science from ``Babe\c{s}-Bolyai'' University, Romania, and the M.S. and Ph.D. degrees in applied mathematics from the University of Maryland, Baltimore County. 
His research interests lie at the intersection of numerical optimization, operations research and high-performance scientific computing with a focus on the optimization of complex energy systems under uncertainty. %He is a member of Society for Industrial and Applied Mathematics (SIAM) and The Institute for Operations Research and the Management Sciences (INFORMS). 
\end{IEEEbiography}
\myvspace{-0.3in}
\begin{IEEEbiography}[{\includegraphics[width=1in,clip,keepaspectratio]{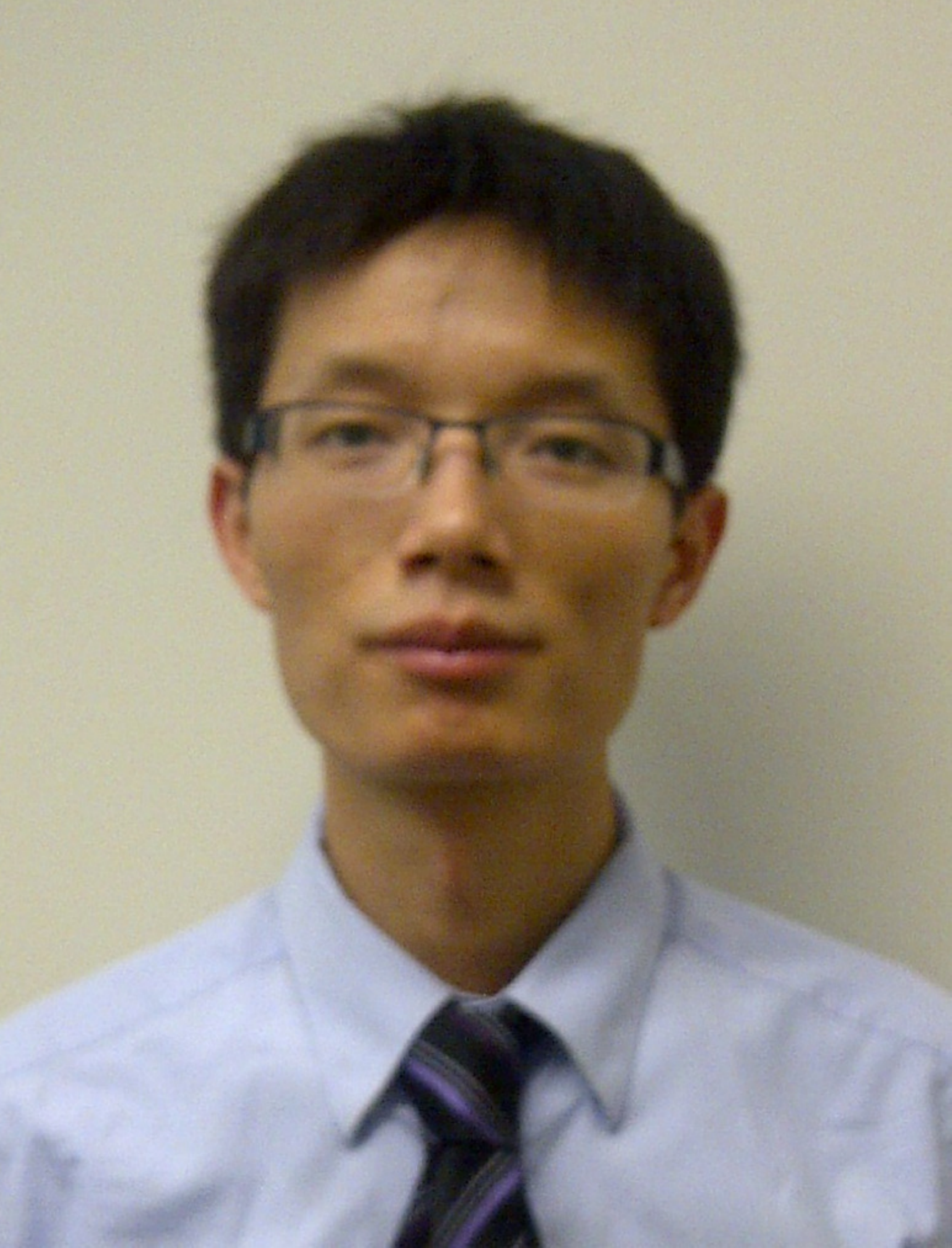}}\myvspace{-0.1in}]{Zheng Zhang} received his Ph.D. degree in electrical engineering and computer science from the Massachusetts Institute of Technology (MIT), Cambridge, MA. He is currently a postdoc associate with the Mathematics and Computer Science Division at
Argonne National Laboratory.  His research interests include
high-dimensional uncertainty quantification and data analysis for
nanoscale devices and systems, energy systems, and biomedical
applications. He received the 2015 Doctoral Dissertation Seminar Award
from the Microsystems Technology Laboratory of MIT, and the 2014 Best Paper Award from IEEE Transactions on CAD of Integrated Circuits and Systems. His industrial research experiences include Coventor Inc. and Maxim-IC. %Some of his research results have been implemented in the commerical MEMS/IC co-design software MEMS+.
\end{IEEEbiography}
\myvspace{-0.3in}
\begin{IEEEbiography}[{\includegraphics[width=1in,clip,keepaspectratio]{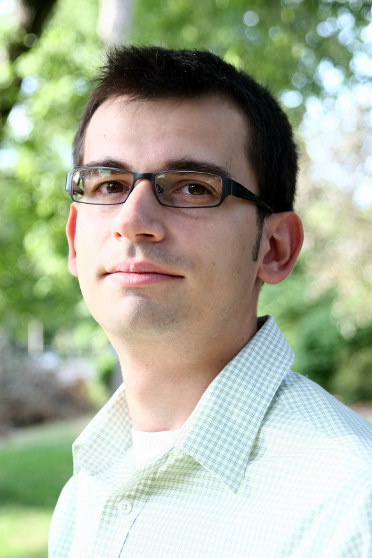}}]{Emil M. Constantinescu} received his Ph.D. degree in computer science
from Virginia Tech, Blacksburg, in 2008. He is currently a computational
mathematician in the Mathematics and Computer Science Division at
Argonne National Laboratory and he is on the editorial board of SIAM Journal on Scientific Computing. His research interests include
numerical analysis of time-stepping algorithms and their applications to
energy systems.\end{IEEEbiography}
\myvspace{-0.3in}
\begin{IEEEbiography}[{\includegraphics[width=1in,height=1.25in,clip,keepaspectratio]{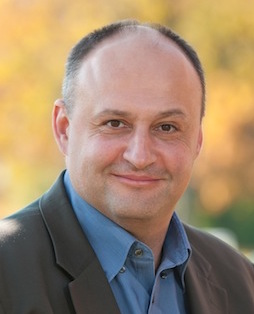}}]{Mihai Anitescu} 
 is a senior computational mathematician in the Mathematics and
 Computer Science Division at Argonne National Laboratory, a professor
 in the Department of Statistics at the University of Chicago, and a
 Senior Fellow of the Computation Institute at the University of
 Chicago. He obtained his engineer diploma (electrical engineering)
 from the Polytechnic University of Bucharest in 1992 and his Ph.D. in
 applied mathematical and computational sciences from the University
 of Iowa in 1997. He specializes in the areas of numerical optimization, computational science, numerical analysis and uncertainty quantification. He co-authored more than 100 peer-reviewed papers in scholarly journals, book chapters, and conference proceedings, and he is on the editorial board of Mathematical Programming A and B, SIAM Journal on Optimization, SIAM Journal on Scientific Computing, and SIAM/ASA Journal in Uncertainty Quantification; and he is a senior editor for Optimization Methods and Software. 
\end{IEEEbiography}
%\myvspace{-0.4in}
%\clearpage
%\newpage
\vfill
\begin{flushright}
\scriptsize \framebox{\parbox{3.2in}{(To be removed before publication)
    Government License: The submitted manuscript has been created by
    UChicago Argonne, LLC, Operator of Argonne National Laboratory
    (``Argonne").  Argonne, a U.S. Department of Energy Office of
    Science laboratory, is operated under Contract
    No. DE-AC02-06CH11357.  The U.S. Government retains for itself,
    and others acting on its behalf, a paid-up nonexclusive,
    irrevocable worldwide license in said article to reproduce,
    prepare derivative works, distribute copies to the public, and
    perform publicly and display publicly, by or on behalf of the
    Government. }} \normalsize
\end{flushright}

\end{document}